\documentclass[12pt]{amsart}

\usepackage{amsmath,amssymb,amsfonts,amssymb,amsthm}

\usepackage{verbatim}
\usepackage[usenames]{color}
\usepackage{hyperref}
\usepackage{url}
\usepackage{tikz,tikz-qtree,cancel}
\usepackage{array,tikz-qtree,cancel}
\usepackage{graphicx}
\usepackage{adjustbox}
\usepackage{amsthm,graphicx,tikz,appendix,tikz-qtree,ifthen,cancel}
\usetikzlibrary{calc,shapes,patterns,positioning}

\newtheorem{thm}{Theorem}[section]

\newtheorem{cor}[thm]{Corollary}

\newtheorem*{thm7.3}{Theorem \ref{thm.MillikenSWP}}

\theoremstyle{remark}

\theoremstyle{definition}
\newtheorem{defn}[thm]{Definition}

\newtheorem{question}[thm]{Question}
\newtheorem{conjecture}[thm]{Conjecture}

\theoremstyle{remark}

\newcommand{\om}{\omega}

\newcommand{\sse}{\subseteq}

\DeclareMathOperator{\Succ}{Succ}

\DeclareMathOperator{\dom}{dom}
\DeclareMathOperator{\height}{ht}

\newcommand{\re}{\restriction}

\newcommand{\bT}{\mathbb{T}}
\newcommand{\bS}{\mathbb{S}}

\newcommand{\A}{\mathrm{A}}

\newcommand{\bsA}{\boldsymbol{A}}
\newcommand{\bsB}{\boldsymbol{B}}
\newcommand{\bsC}{\boldsymbol{C}}
\newcommand{\bsD}{\boldsymbol{D}}

\newcommand{\ra}{\rightarrow}

\newcommand{\lgl}{\langle}
\newcommand{\rgl}{\rangle}

\newcommand{\Nesetril}{Ne{\v{s}}et{\v{r}}il}

\newcommand{\Rodl}{R{\"{o}}dl}
\newcommand{\Erdos}{Erd{\H{o}}s}
\newcommand{\Posa}{P{\'{o}}sa}

\newcommand{\Fraisse}{Fra{\"{i}}ss{\'{e}}}
\newcommand{\Lauchli}{L{\"{a}}uchli}

\newcommand{\Masulovic}{Ma{\v{s}}ulovi{\'{c}}}
\newcommand{\The}{Th{\'{e}}}

\newcommand{\noprint}[1]{\relax}

\begin{document}

\title[Ramsey Theory on Infinite Structures]{Ramsey Theory on Infinite Structures and the Method of Strong Coding Trees}

\author{Natasha Dobrinen}
\address{University of Denver\\
Department of Mathematics, 2390 S. York St., Denver, CO 80208, USA}
\email{natasha.dobrinen@du.edu}
  \urladdr{\url{http://cs.du.edu/~ndobrine}}

 \thanks{This work was supported by  National Science Foundation Grant DMS-1600781}


\maketitle

\begin{abstract}
This article discusses some recent  trends   in Ramsey theory on infinite structures.
Trees and their Ramsey theory have been vital  to these investigations.
The main ideas behind
 the  author's
recent method of trees with coding nodes
are presented, showing how they can be useful both for coding structures with forbidden configurations as well as those with none.
Using forcing as a tool for finite searches has allowed the
development of
 Ramsey theory on such trees, leading to solutions  for finite big Ramsey degrees of Henson graphs as well as infinite dimensional Ramsey theory of copies of the Rado graph.
Possible future directions for applications of these methods are discussed.
\end{abstract}


\section{Introduction}

Logic and Ramsey theory  have  a long  interconnected history.
In  1929, Ramsey proved his
 celebrated theorem in order to obtain a partial solution to
Hilbert's Entscheidungsproblem.
This problem, posed by Hilbert in 1928,  asked for
an algorithm which decides the validity of  any statement of first-order logic.
By G\"{o}del's
Completeness Theorem (1929), this is equivalent to asking for an algorithm which can decide whether a given formula is provable from the logical axioms.
Ramsey applied his partition theorem on the natural numbers  to show that the validity of
  formulas with only universal quantifiers in  normal form   is decidable \cite{Ramsey30}.
In light of the  Church-Turing thesis
(\cite{Church36} and \cite{Turing36}) which precludes any  complete solution  to Hilbert's Entscheidungsproblem,
Ramsey's result is all the more striking in its success.

Over the decades, this interplay between Ramsey theory and logic has continued, each subject motivating, informing, and providing techniques to solve problems in the other.
An instance of this is seen in the 1966 result of Halpern and \Lauchli\ \cite{Halpern/Lauchli66}.
While  investigating  the problem of whether or not the Boolean Prime Ideal Theorem (BPI) is strictly weaker than the Axiom of Choice (AC) over the Zermelo-Fraenkel Axioms (ZF),
Halpern and \Lauchli\  proved a theorem which was later interpreted to be a Ramsey theorem on products of finitely many  trees.
This theorem was  central to  the proof of Halpern and
L\'{e}vy  that, indeed, BPI is strictly weaker than AC over ZF \cite{Halpern/Levy71}.
In turn, the Halpern-\Lauchli\ theorem  provided means for proving Ramsey theorems for colorings of finite structures inside of infinite structures, such as the rationals as an ordered structure  and the Rado graph.

Harrington later produced a  proof of the Halpern-\Lauchli\ Theorem using
 the set-theoretic method of forcing.
His approach was novel  in  that the language and techniques of forcing are used  to do unbounded searches for finite objects.
This is in contrast to the more common use of forcing to
 obtain ZFC results by proving
the
 existence  of a $\Sigma^1_2$ definable object  in a generic extension, and then applying Shoenfield absoluteness to deduce that the object must be in the ground model.
In an interesting turn of events, Harrington's method of proof provided the backdrop for recent
work of the author in \cite{DobrinenJML20}, \cite{DobrinenH_k19}, and \cite{DobrinenRado19}.
These results will be discussed in Section \ref{sec.bRd}
and
ideas of key methods developed  to obtain these results will be set forth in Section \ref{sec.sct}.

This article discusses some recent  trends   and future directions in Ramsey theory on infinite structures.
 Section \ref{sec.overview}  recalls  Ramsey theory  on the natural numbers, and then
 reviews how these ideas extend to  structures.
 Ties between Ramsey theory on structures and
 topological dynamics, due to Kechris, Pestov, and Todorcevic in \cite{Kechris/Pestov/Todorcevic05} and a recent development due to Zucker in \cite{Zucker19}
 provide additional motivation for these investigations.
 In Section
\ref{sec.bRd}, we present an overview of   big Ramsey degrees and infinite dimensional Ramsey theory on infinite structures.
As the focus here  is Ramsey theory on infinite structures,  and as the literature on finite structures is vast, we do not even attempt to do justice to Ramsey theory of finite structures in this article.
Only a few of those results will be stated in order to provide the reader with some intuition.

The main tools
used so far  to obtain Ramsey results  on infinite structures, aside from Ramsey's Theorem itself,
have been
 Ramsey theorems on trees
in the vein of  Halpern and \Lauchli,
and very recently, category theory (see \cite{Masulovic18}).
In
 Section
\ref{sec.sct}, we present
 an overview of the  recent method of trees with coding nodes, first developed by the author in \cite{DobrinenJML20} to prove finite big Ramsey degrees for the universal ultrahomogeneous triangle-free graph.
 This method was extended in
\cite{DobrinenH_k19} and
 \cite{DobrinenRado19}
 to determine Ramsey theory on all the Henson graphs as well as infinite dimensional Ramsey theory on the Rado graph.
 Forcing seemed the best first approach in the  search for Ramsey theorems on these trees, and it turned out to work.
An overview of some of the ideas involved in these proofs are provided.
In Section
\ref{sec.fd} we
 point to future directions in Ramsey theory of infinite structures where
  these methods are
  likely to prove efficacious.


\section{Finite and infinite dimensional   Ramsey theory:\\ from  natural numbers to structures}\label{sec.overview}

Ramsey's theorem for the natural numbers is the following:

\begin{thm}[Infinite Ramsey Theorem, \cite{Ramsey30}]\label{thm.RamseyInfinite}
Given integers $k,r\ge 1$ and a coloring $c:[\om]^k\ra r$,
 there is an infinite set  $N\sse \om$
such that $c\re [N]^k$ is constant.
\end{thm}

Here and throughout, we use the set-theoretic notation $[X]^k$ to denote the collection of all subsets of $X$ with exactly $k$ members.
Using the Hungarian arrow notation, Theorem \ref{thm.RamseyInfinite} is written as
\begin{equation}
\om\ra(\om)^k_r.
\end{equation}
 The set $N$ is said to be {\em monochromatic}  or  {\em homogeneous}.
Using compactness, one obtains the finite version of Ramsey's Theorem.

\begin{cor}[Finite Ramsey Theorem, \cite{Ramsey30}]\label{cor.RamseyFinite}
Given integers $k,r,m\ge 1$  with
$k\le m$, there is an integer $n>m$ such that for any
coloring $c:[n]^k\ra r$,
there is a subset $N\sse n$ of cardinality $m$ such that $c\re [N]^k$ is constant.
\end{cor}

Theorem \ref{thm.RamseyInfinite} and Corollary \ref{cor.RamseyFinite}
 are referred to as  {\em  finite dimensional Ramsey theory}, since the
sets being colored in these theorems all have the same  fixed finite size.
{\em Infinite dimensional Ramsey theory} is concerned with coloring infinite sets of natural numbers.
Although the Axiom of Choice implies  there is a  coloring of all infinite subsets of natural numbers into two colors for which there is no infinite monochromatic set (see \Erdos-Rado \cite{Erdos/Rado52}),
for any coloring  which induces a sufficiently definable partition of the Baire space, monochromatic subsets are to be found.
In the context of Ramsey theory, the simplest representation of the Baire space is
 $[\om]^{\om}$,
the set of all infinite subsets of $\om$,  endowed
with the Tychonoff topology.
Nash-Williams  initiated the study of  infinite dimensional Ramsey theory in 1965,
proving that for  any clopen  set $\mathcal{C}$ in the Baire space, there is some infinite set $X$ such that
 $[X]^{\om}$ is either contained in $\mathcal{C}$ or else is disjoint from it \cite{NashWilliams65}.

A few years later, Galvin and Prikry  extended this to Borel sets in a strong way.
The following notation is central to studies of infinite dimensional Ramsey theory:  For $s\in[\om]^{\om}$ and $A\in[\om]^{\om}$, define
 \begin{equation}
 [s,A]=\{B\in[\om]^{\om}:s\sqsubset B\sse A\},
 \end{equation}
where $s$ is finite and $s\sqsubset B$ means that $s$ is an initial segment of $B$, given their strictly increasing enumerations.
We use the terminology of \cite{TodorcevicBK10}.

\begin{defn}\label{defn.Ramsey}

A subset $\mathcal{X}\sse[\om]^{\om}$ is {\em Ramsey} if for each  non-empty set $[s,A]$,
there is a $B\in [s,A]$ such that either $[s,B]\sse\mathcal{X}$ or else
$[s,B]\cap\mathcal{X}=\emptyset$.
\end{defn}

\begin{thm}[Galvin and Prikry, \cite{Galvin/Prikry73}]\label{thm.GP}
Every  Borel subset of the  Baire space is Ramsey.
\end{thm}

Soon after this, Silver  showed that analytic sets are Ramsey \cite{Silver71}.
Mathias \cite{Mathias77} and Louveau \cite{Louveau74}   attained infinite dimensional Ramsey  results for the case where the infinite set comes from a Ramsey ultrafilter, a related and rich area which is not the focus of  this article.
  The pinnacle of infinite dimensional Ramsey theory on the Baire space  was achieved by Ellentuck in 1974, who provided a topological characterization of those sets which are Ramsey.
 The {\em Ellentuck topology} is the topology on $[\om]^{\om}$ generated by basic open
sets of the form
$ [s,A]$.
This topology refines the Tychonoff topology on the Baire space.

\begin{thm}[Ellentuck, \cite{Ellentuck74}]\label{thm.Ellentuck}
A subset $\mathcal{X}\sse [\om]^{\om}$ is Ramsey if and only if it has the property of Baire with respect to the Ellentuck topology.
\end{thm}

Such theorems can be notated by
\begin{equation}
\om\ra_*(\om)^{\om}_2,
\end{equation}
where $\ra_*$ denotes that the sets are partitioning $[\om]^{\om}$ are required to be  definable in some sense  with respect to some topology (see Section 11 of \cite{Kechris/Pestov/Todorcevic05}).

Overlapping the development of infinite dimensional Ramsey theory on the Baire space,
Ramsey theory  on relational structures  began to unfold.
We review only the basics of \Fraisse\ theory for  relational structures;
more general background  on \Fraisse\ theory can be found in \Fraisse's original paper  \cite{Fraisse54},   as well as  \cite{Kechris/Pestov/Todorcevic05}.
We call $\mathcal{L}=\{R_i\}_{i\in I}$  a
{\em  relational signature} if is  a (countable) collection of  {\em relation symbols} $R_i$, where $n(i)$
 denotes the {\em arity} of $R_i$, for each $i\in I$.
A {\em structure} for $\mathcal{L}$
is a structure
 $\bsA=\lgl |\bsA|, \{R_i^{\bsA}\}_{i\in I}\rgl$,
where $|\bsA|\ne\emptyset$ is the {\em universe} of $\bsA$  and for each $i\in I$, $R_i^{\bsA}\sse A^{n(i)}$.
An {\em embedding} between structures $\bsA,\bsB$ for $\mathcal{L}$ is an injection $\iota:|\bsA|\ra|\bsB|$ such that
for all $i\in I$, $R_i^{|\bsA|}(a_1,\dots,a_{n(i)})\leftrightarrow R_i^{|\bsB|}(\iota(a_1),\dots,\iota(a_{n(i)}))$.
If $\iota$ is the identity map, then we say that $\bsA$ is a {\em substructure} of $\bsB$.
An  {\em isomorphism}  is an  embedding which is onto it image.
We write $\bsA\le\bsB$  exactly when  there is an embedding of  $\bsA$  into $\bsB$;
 $\bsA\cong\bsB$  denotes that $\bsA$ and $\bsB$ are isomorphic.

A class $\mathcal{K}$ of finite structures for a relational signature $\mathcal{L}$
 is called a {\em \Fraisse\ class}  if it is hereditary, satisfies the joint embedding and amalgamation properties, contains (up to isomorphism) only countably many structures, and contains structures of arbitrarily large finite cardinality.
These notions are recalled here for the reader's convenience.
 $\mathcal{K}$ is  {\em hereditary} if whenever $\bsB\in\mathcal{K}$ and  $\bsA\le\bsB$, then also $\bsA\in\mathcal{K}$.
$\mathcal{K}$ satisfies the {\em joint embedding property} if for any $\bsA,\bsB\in\mathcal{K}$,
there is a $\bsC\in\mathcal{K}$ such that $\bsA\le\bsC$ and $\bsB\le\bsC$.
 $\mathcal{K}$ satisfies the {\em amalgamation property} if for any embeddings
$f:\bsA\ra\bsB$ and $g:\bsA\ra\bsC$, with $\bsA,\bsB,\bsC\in\mathcal{K}$,
there is a $\bsD\in\mathcal{K}$ and  there are embeddings $r:\bsB\ra\bsD$ and $s:\bsC\ra\bsD$ such that
$r\circ f = s\circ g$.

Let $\mathcal{K}$  be a \Fraisse\ class.
For $\bsA,\bsB\in\mathcal{K}$ with $\bsA\le\bsB$, we use ${\bsB\choose\bsA}$ to denote the set of all substructures of $\bsB$ which are isomorphic to $\bsA$.
Given structures $\bsA\le\bsB\le\bsC$ in $\mathcal{K}$, we write
$$
\bsC\ra(\bsB)_k^{\bsA}
$$
to denote that for each coloring of ${\bsC\choose \bsA}$ into $k$ colors, there is a $\bsB' \in {\bsC\choose\bsB}$ such that
${\bsB'\choose\bsA}$ is  {\em monochromatic}, meaning
 every member of ${\bsB'\choose\bsA}$ has the same color.

 \begin{defn}\label{defn.RP}
A \Fraisse\ class  $\mathcal{K}$ has the {\em Ramsey property} if  for any two structures $\bsA\le\bsB$ in $\mathcal{K}$ and any  $k\ge 2$,
there is a $\bsC\in\mathcal{K}$ with $\bsB\le\bsC$ such that
$\bsC\ra (\bsB)^{\bsA}_k$.
\end{defn}

Investigations into which \Fraisse\ classes have the Ramsey property commenced when
Graham and Rothschild proved in 1971 that the class of finite Boolean algebras has the Ramsey property \cite{Graham/Rothschild71}.
Soon after this,
Graham, Leeb, and Rothschild showed that  the class of
finite vector spaces over a finite field  have the Ramsey property \cite{Graham/Leeb/Rothschild72},
\cite{Graham/Leeb/Rothschild73}.
Several years later,
the class of
 finite  ordered graphs were  found to have the
  Ramsey property;
  this was proved independently by
  Abramson and Harrington \cite{Abramson/Harringon78} and  by \Nesetril\ and \Rodl\ \cite{Nesetril/Rodl77}, \cite{Nesetril/Rodl83}.
The main theorem of the papers \cite{Nesetril/Rodl77} and \cite{Nesetril/Rodl83} furthermore  proved
 that  all set-systems of  finite ordered relational structures omitting some irreducible substructure have the Ramsey property.
This includes
 the  classes of finite ordered graphs omitting $k$-cliques, denoted $\mathcal{G}_k^{<}$, for each $k\ge 3$.
Over the past several decades,  more \Fraisse\ classes were  shown to have the Ramsey property.
The correspondence between the Ramsey property and
extreme amenability, proved by Kechris, Pestov, and Todorcevic in 2005 in \cite{Kechris/Pestov/Todorcevic05},
has propelled  a recent
burst of discoveries of more \Fraisse\ classes with the Ramsey property.

It is interesting that while most \Fraisse\   classes of finite unordered structures do not have
the Ramsey property, often
 equipping
the class with an additional linear order  produces  the Ramsey property.
In such cases,  some
 remnant of the Ramsey property  remains in the unordered reduct.
Following notation  in \cite{Kechris/Pestov/Todorcevic05},
given a \Fraisse\ class $\mathcal{K}$,
for  $\bsA\in\mathcal{K}$,
$t(\bsA,\mathcal{K})$ denotes  the smallest number $t$, if it exists, such that
for each $\bsB\in \mathcal{K}$ with $\bsA\le \bsB$ and for each $j\ge 2$,
there is some $\bsC\in\mathcal{K}$ into which $\bsB$ embeds  such that for
 any coloring $c:{\bsC \choose \bsA}\ra j$,
 there is a $\bsB'\in {\bsC \choose \bsB}$ such that the restriction of $c$ to ${\bsB'\choose \bsA}$ takes no more than $t$ colors.
In the arrow notation, this is written as
\begin{equation}
\bsC\ra (\bsB)^{\bsA}_{j,t(\bsA,\mathcal{K})}.
\end{equation}
A class  $\mathcal{K}$ has  {\em finite (small) Ramsey degrees} if
for each $\A\in\mathcal{K}$ the number
 $t(\bsA,\mathcal{K})$  exists.
The number $t(\bsA,\mathcal{K})$ is called the {\em Ramsey degree of $\bsA$} in $\mathcal{K}$ \cite{Fouche98}.
Notice that $\mathcal{K}$ has the Ramsey property if and only if $t(\bsA,\mathcal{K})=1$ for each $\bsA\in\mathcal{K}$.
The  connection between \Fraisse\ classes with finite Ramsey degrees and ordered expansions
was initiated in 
 Section 10 of \cite{Kechris/Pestov/Todorcevic05}
 and extended to precompact  expansions in  \cite{NVTHabil}.
The existence of a Ramsey expansions 
 were later
shown to be equivalent to small Ramsey degrees in  \cite{Zucker16}.
In particular, the \Fraisse\
 classes of  finite  unordered graphs and finite unordered graphs omitting $k$-cliques have   finite small Ramsey degrees.

At this point, it is natural to ask the following:
Which infinite structures have  properties similar  Theorem \ref{thm.RamseyInfinite}?
As  it is often not possible to obtain one color for all copies of a given object,  the following definition extends the notion of small Ramsey degrees (rather than Ramsey property) to infinite structures.

\begin{defn}[\cite{Kechris/Pestov/Todorcevic05}]\label{defn.bRd}
Given an infinite structure $\mathcal{S}$ and a finite substructure $\bsA\le \mathcal{S}$,
let $T(\bsA,\mathcal{S})$ denote the least integer $T\ge 1$, if it exists, such that
given any coloring of ${\mathcal{S}\choose \bsA}$ into finitely many colors, there is a
 substructure $\mathcal{S}'$ of $\mathcal{S}$, isomorphic to $\mathcal{S}$,  such that ${\mathcal{S}'\choose \bsA}$ takes no more than $T$ colors.
This may be written succinctly as
\begin{equation}\label{eq.bRd}
\forall j\ge 1,\ \ {\mathcal{S}}\ra ({\mathcal{S}})^{\bsA}_{j,T(\bsA,\mathcal{S})}.
\end{equation}
We say that
 $\mathcal{S}$ has {\em finite big Ramsey degrees} if for each finite substructure
$\bsA\le\mathcal{S}$,
there is an integer $T(\bsA,\mathcal{S})$
such that (\ref{eq.bRd}) holds.
\end{defn}

By a theorem of  Hjorth in \cite{Hjorth08},  if $\mathcal{K}$ is a \Fraisse\ class with 
limit $\mathbb{F}$ such that $|$Aut$(\mathbb{F})|>1$,
then there is some $\bsA\in\mathcal{K}$ whose big Ramsey degree is larger than one (possibly not finite). 
Thus any non-rigid \Fraisse\ structure cannot have an exact analogue of the Infinite Ramsey's Theorem, where all big Ramsey degrees are one.

Some authors prefer to  color 
 embeddings   of a given  $\bsA$ into $\mathbb{F}$.  
The relationship between them is  simple:
A structure  $\bsA\in\mathcal{K}$ has (small or big) Ramsey degree $\ell$ for copies if and only if $\bsA$ has (small or big) Ramsey degree $\ell\cdot |$Aut$(\bsA)|$ for embeddings. 

The following  question  has been investigated for several decades.

\begin{question}\label{q.bRd}
Which infinite structures have finite big Ramsey degrees?
\end{question}

This question has been asked most often in regard to
 \Fraisse\  limits of \Fraisse\ classes of finite structures,
 which we shall call simply {\em \Fraisse\ structures}.
  A
\Fraisse\ structure $\mathbb{F}$ for a \Fraisse\ class $\mathcal{K}$
is a  countably infinite structure which is {\em universal} for $\mathcal{K}$ (each member of $\mathcal{K}$ embeds into $\mathbb{F}$)
and {\em ultrahomogeneous} (any isomorphism between two finite substructures of $\mathbb{F}$ extends to an automorphism of $\mathbb{F}$).
However, Question \ref{q.bRd} also makes sense in the context of infinite structures which are universal for some $\mathcal{K}$ but not ultrahomogeneous, and we will see some recent progress in this direction as well in the next section.

While Question  \ref{q.bRd} has been
of interest for many decades,
it  gained renewed traction  in the early 2000's, both
because of
results on the big Ramsey degree of the Rado graph in
 \cite{Sauer06} and \cite{Laflamme/Sauer/Vuksanovic06},
and
 because of  Kechris, Pestov, and Todorcevic's   question about  finding  an analogue of their correspondence for \Fraisse\ structures which have finite big Ramsey degrees (see  Section 11 of  \cite{Kechris/Pestov/Todorcevic05}).
Such an analogue was obtained by Zucker in \cite{Zucker19}.
Section \ref{sec.bRd} will give an overview of known results on big Ramsey degrees of infinite structures.

Another natural  question is the following:
Which infinite structures carry
analogues of Theorems \ref{thm.GP} or \ref{thm.Ellentuck}?

\begin{question}\label{question.2.9}
For which infinite structures $\mathcal{S}$ does
\begin{equation}
\mathcal{S}\ra_*(\mathcal{S})^{\mathcal{S}}_2
\end{equation}
hold, where $\ra_*$ denotes that the partition is suitably definable, given some natural topology on the space ${\mathcal{S}\choose\mathcal{S}}$?
\end{question}

This question was brought to light in Problem 11.2 of \cite{Kechris/Pestov/Todorcevic05}, which  asked for an analogue of the KPT correspondence  for \Fraisse\ structures with some infinite dimensional Ramsey theory,
while simultaneously pointing out that very little was known in this direction.
In the  recent paper \cite{DobrinenRado19},
the author showed that for the collection of subcopies of the Rado
graph with a certain tree-structural property,  all Borel sets are Ramsey.
This will be  discussed in  Section \ref{sec.sct},  where  trees with coding nodes will be  presented.


\section{Big Ramsey degrees   on infinite structures:\\  An overview of previous results and
methods}\label{sec.bRd}

Ramsey theory on infinite structures seeks to find  out which infinite structures carry analogues of Ramsey's Theorem.
The first
line of inquiry  investigates
which infinite structures are
 {\em indivisible}, meaning that given any partition of the universe of the structure into finitely many pieces, one of the pieces contains a copy of the structure.
 This is the same as saying
 (in the terminology  of Definition \ref{defn.bRd})
  that
 substructures of size one have
 finite  big Ramsey
  degree one.

It is straightforward to see that given a
partition of the rationals into finitely many pieces, one of the pieces contains a copy of the rationals, that is, a dense linear order (without endpoints).
Thus, the rationals as a linearly ordered structure is indivisible.
The
  Rado graph, or random graph,  is also indivisible.
  This was proved by Henson in \cite{Henson71},
using  the  extension property of the Rado graph.

 In contrast, proving the indivisibility of the generic $k$-clique-free graphs, denoted $\mathcal{H}_k$, required ideas beyond  their extension properties.
These {\em Henson graphs} were first constructed by
 Henson in 1971  in \cite{Henson71} as subgraphs of the Rado graph.
 In hindsight, these  are seen to be
 \Fraisse\ limits of the \Fraisse\ classes
 of  finite $k$-clique-free graphs.
In \cite{Henson71},
Henson proved  that $\mathcal{H}_k$ is {\em weakly indivisible}:
Given a partition of the vertices of $\mathcal{H}_k$ into two pieces, either the first piece of the partition contains a copy of $\mathcal{H}_k$, or else the second piece contains a copy of each finite $k$-clique-free graph.
It is interesting to note that this was proved several years before \Nesetril\ and \Rodl's result  in \cite{Nesetril/Rodl77} that the  \Fraisse\ classes of ordered $k$-clique-free graphs have the Ramsey property.
The first result on the indivisibility of Henson graphs was due to
Komj\'{a}th and \Rodl, in  \cite{Komjath/Rodl86}, where  they proved  that  $\mathcal{H}_3$  is indivisible.
Soon after,  El-Zahar and Sauer extended this result to all Henson graphs
 in \cite{El-Zahar/Sauer89}.

 A sample of other notable results  regarding indivisibility  are the following:
 Hindman showed, as
a consequence of his partition theorem for finite sums of natural numbers, that
 the  vector space  of countable dimension over the finite field $\mathbb{F}_2$ is
  indivisible  \cite{Hindman74}.
  However, a  vector space  of countable dimension over any other (nontrivial)  finite field is not.
  A proof of this folklore theorem appears in   \cite{Laflamme/NVT/Sauer/Pouzet11}, where 
   Laflamme, Nguyen Van Th\'{e}, Pouzet, and  Sauer
  prove that all  such  vector spaces  retain a  slightly weaker property called weak indivisibility. 
 Nguyen Van \The\  and Sauer proved  that for each $m\ge 1$,  the countable Urysohn metric  space with distances in $\{1,\dots,m\}$ is indivisible \cite{NVTSauer09}.
 Later, Sauer showed that for any finite set of distances, the Urysohn space with that distance set is indivisible \cite{Sauer12}.
For more on indivisible structures, the reader is referred to the excellent Habilitation of Nguyen Van \The\
 \cite{NVTHabil}.

While the study of indivisibility of infinite structures continues to be a rich and challenging subject,
it is interesting that
the nascence of the broader subject of big Ramsey degrees  can actually be traced back to an example of
  Sierpi\'{n}ski.
 Considering the  rationals as a linearly  ordered structure,
 Sierpi\'{n}ski
  showed that there is a coloring of pairs of rationals into two colors such that any infinite subset which is again a dense linear order preserves both colors on its pairsets.
 His coloring plays the linear ordering of the rationals against a well-ordering of the rationals as follows, and colors pairs  of rationals red if the two orders agree on the pair, and blue otherwise.
Both colors persist in every subset of the rationals forming a dense linear order.

 This phenomenon is perhaps best seen via  trees.
  For $s,t\in 2^{<\om}$, define  $s \triangleleft t$  iff one of the following holds:
(a)
$s<_{\mathrm{lex}} t$,
(b)
$s\sqsubset t$ and $t(|s|)=1$, or
(c)
$t\sqsubset s$ and $s(|t|)=0$.
Notice that $(2^{<\om},\triangleleft)$ is isomorphic to $(\mathbb{Q},<)$.
  Define the coloring
  \begin{equation}
c(\{s,t\})=
\begin{cases}
0 & \text{if\ } |s|\le |t|\mathrm{\ and\ } s\triangleleft t\\
1 & \text{otherwise}
\end{cases}
\end{equation}
Then given any subset $S\sse 2^{<\om}$ for which $(S,\triangleleft)\cong (\mathbb{Q},<)$,
both colors will persist in $S$; hence $T(2,\mathbb{Q})\ge 2$.

In his  PhD thesis \cite{DevlinThesis}, Devlin  found the precise big Ramsey degrees for finite sets of rationals, building on prior unpublished work of Galvin, who proved that $T(2,\mathbb{Q})=2$, and of Laver, who proved existence of $T(k,\mathbb{Q})$ for all natural numbers $k$.
Here, we will not  provide much background or history, as there are already
thorough expositions of Devlin's results
 in \cite{Farah/TodorcevicBK} and
\cite{TodorcevicBK10}.

A key component in  finding upper bounds for the numbers $T(k,\mathbb{Q})$ is a Ramsey theorem for trees, due to Milliken.
This actually turns out to be at the heart of many big Ramsey degree results.
Milliken's theorem utilizes the following
notion of strong tree and a theorem of  Halpern and \Lauchli\ which we now briefly review.
A historical record of the Halpern-\Lauchli\ Theorem and its variants can be found in  \cite{Larson12}.
For $t\in \om^{<\om}$,   $|t|=\dom(t)$ is  the  {\em  length of $t$}.
In the subject of Ramsey theory on trees,
we say that a subset $T\sse \om^{<\om}$ is  a {\em tree} if
 there is a subset $ L\sse\om$ such that $T=\{t\re l: t\in T,\ l\in L\}$.
 Thus, in this definition, a tree is closed under initial segments at levels of the tree, but it is not necessarily  closed under all initial segments in $\om^{<\om}$.
For $t\in T$,  the {\em  height of $t$}, $\height_T(t)$, is the order-type of the set  $\{u\in T:u\subset t\}$, linearly ordered by $\sse$.
We write
 $T(n)$ to  denote $\{t\in T:\height_T(t)=n\}$.
 For $t\in T$,  let $\Succ_T(t)=\{u\re(|t|+1):u\in T$ and $u\supset t\}$, noting that $\Succ_T(t)$ will not be  a set of nodes in $T$  if $T$ does not contain any nodes of length $|t|+1$.
 For $s,t\in T$ with $|s|<|t|$, the number $t(|s|)$ is called the {\em passing number} of $t$ at $s$.
 Passing numbers are used to code information about the binary relations
 satisfied by 
 the members of (the universe of) a given structure represented by $s$ and $t$.

\begin{defn}
Let $T\sse \om^{<\om}$ be a  finitely branching tree.
A subset $S\sse T$ is a
{\em  strong subtree of $T$}  if and only if
there is  an increasing sequence of natural numbers
$\lgl m_n:n<N\rgl$ ,  where $N\le \om$, such
that
$S=\bigcup_{n<N}S(n)$, and
 for each $n<N$,
\begin{enumerate}
\item
 $S(n)\sse T(m_n)$, and
\item
$s\in S(n)$ and $u\in\Succ_T(s)$,
there is exactly one node in  $S(n+1)$ extending $u$.
\end{enumerate}
Given $k\ge 1$, we say that $S$ is a {\em $k$-strong subtree}  of $T$ if $N=k<\om$.
\end{defn}

For the next theorem, define the notation:
\begin{equation}
\bigotimes_{i<d}T_i:=\bigcup_{n<\om}  \prod_{i<d} T_i(n).
\end{equation}

The following is the strong tree version of the Halpern-\Lauchli\ Theorem.

\begin{thm}[Halpern-\Lauchli, \cite{Halpern/Lauchli66}]\label{thm.HL}
Let  $T_i\sse \om^{<\om}$, $i<d$, be  finitely branching trees
 with no terminal nodes and  let  $r\ge 2$.
Given a coloring
$c:\bigotimes_{i<d} T_i \ra r$,
there is an increasing sequence $\lgl m_n:n<\om\rgl$ and
 strong  subtrees $S_i\le T_i$
 such that for all $i<d$ and $n<\om$, $S_i(n)\sse T_i(m_n)$, and
$c$ is constant on $\bigotimes_{i<d}S_i$.
\end{thm}

The following theorem  is
proved by inductively applying Theorem \ref{thm.HL} $\om$-many times.

\begin{thm}[Milliken, \cite{Milliken79}]\label{thm.Milliken}
Let $T\sse \om^{<\om}$ be a  strong tree with no maximal nodes.
Let  $k\ge 1$, $r\ge 2$, and $c$ be a coloring of all
$k$-strong subtrees of
 $T$ into $r$ colors.
Then there is a  strong subtree $S\sse T$ such that all  $k$-strong subtrees of  $S$ have the same color.
\end{thm}


\begin{figure}
\begin{tikzpicture}[grow'=up,scale=.3]
\tikzstyle{level 1}=[sibling distance=4in]
\tikzstyle{level 2}=[sibling distance=2in]
\tikzstyle{level 3}=[sibling distance=1in]
\node {} coordinate (t9)
child{coordinate (t0) edge from parent[thick, color=red]
			}
		child{ coordinate(t1) edge from parent[thick, color=red]
			child{ coordinate(t10)
			child{ coordinate(t100)edge from parent[draw=none]
}
child{  coordinate(t101)
}
}
			child{  coordinate(t11)edge from parent[draw=none]
} };

\node[circle, fill=red,inner sep=0pt, minimum size=5pt] at (t0) {$s$};
\node[circle, fill=red,inner sep=0pt, minimum size=5pt] at (t101) {$t$};

\end{tikzpicture}
\hspace{20pt}
\begin{tikzpicture}[grow'=up,scale=.3]
\tikzstyle{level 1}=[sibling distance=4in]
\tikzstyle{level 2}=[sibling distance=2in]
\tikzstyle{level 3}=[sibling distance=1in]
\node {} coordinate (t9)
child{coordinate (t0) edge from parent[thick, color=blue]
			child{coordinate (t00) edge from parent[color=blue]
			child{ coordinate(t000)edge from parent[draw=none]}
child{ coordinate(t001)  edge from parent[color=blue]
}}
			child{ coordinate(t01)edge from parent[draw=none]
}}
		child{ coordinate(t1) edge from parent[thick, color=blue] };

\node[circle, fill=blue,inner sep=0pt, minimum size=5pt] at (t001) {$s$};
\node[circle, fill=blue,inner sep=0pt, minimum size=5pt] at (t1) {$t$};

\end{tikzpicture}
\caption{The two similarity types of pairs of rationals}
\end{figure}
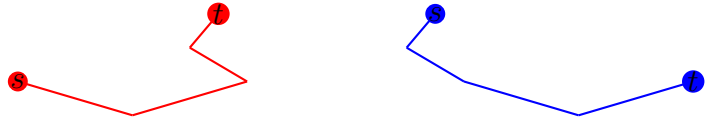


\begin{figure}
\begin{tikzpicture}[grow'=up,scale=.3]
\tikzstyle{level 1}=[sibling distance=4in]
\tikzstyle{level 2}=[sibling distance=2in]
\tikzstyle{level 3}=[sibling distance=1in]
\node {} coordinate (t9)
child{coordinate (t0) edge from parent[thick, color=red]
			child{coordinate (t00)edge from parent[color=black]
child{coordinate (t000)
}
child{ coordinate(t001) edge from parent[draw=none]
}}
			child{ coordinate(t01)edge from parent[color=black]
child{ coordinate(t010) edge from parent[draw=none]
}
child{ coordinate(t011) edge from parent[color=black]
}}}
		child{ coordinate(t1) edge from parent[thick, color=red]
			child{ coordinate(t10)
child{ coordinate(t100) edge from parent[draw=none]
}
child{ coordinate(t101)
}}
			child{  coordinate(t11)edge from parent[color=black]
child{ coordinate(t110)edge from parent[draw=none]
}
child{  coordinate(t111)
}} };

\node[circle, fill=red,inner sep=0pt, minimum size=7pt] at (t0) {$s$};
\node[circle, fill=red,inner sep=0pt, minimum size=7pt] at (t101) {$t$};

\node[circle, draw, inner sep=0pt, minimum size=7pt] at (t9) {};

\node[circle, draw, inner sep=0pt, minimum size=7pt] at (t000) {};

\node[circle, draw, inner sep=0pt, minimum size=7pt] at (t011) {};

\node[circle, draw, inner sep=0pt, minimum size=7pt] at (t1) {};

\node[circle, draw, inner sep=0pt, minimum size=7pt] at (t111) {};

\end{tikzpicture}
\hspace{20pt}
\begin{tikzpicture}[grow'=up,scale=.3]
\tikzstyle{level 1}=[sibling distance=4in]
\tikzstyle{level 2}=[sibling distance=2in]
\tikzstyle{level 3}=[sibling distance=1in]
\node {} coordinate (t9)
child{coordinate (t0) edge from parent[thick, color=blue]
			child{coordinate (t00) edge from parent[draw=blue]
child{coordinate (t000) edge from parent[draw=none]}
child{coordinate(t001) edge from parent[color=blue]}}			
			child{ coordinate(t01)edge from parent[color=black]
child{ coordinate(t010)}
child{ coordinate(t011) edge from parent[draw=none]
}}}
		child{ coordinate(t1) edge from parent[thick, color=blue]
child{coordinate (t10)edge from parent[draw=black]
child{coordinate(t100)}
child{coordinate(t101)edge from parent[draw=none]}}
child{coordinate (t11)edge from parent[draw=black]
child{coordinate (t110)edge from parent[draw=none]}
child{coordinate(t111)}}		
};

\node[circle, fill=blue,inner sep=0pt, minimum size=5pt] at (t001) {$s$};
\node[circle, fill=blue,inner sep=0pt, minimum size=5pt] at (t1) {$t$};

\node[circle, draw, inner sep=0pt, minimum size=7pt] at (t9) {};
\node[circle, draw, inner sep=0pt, minimum size=7pt] at (t0) {};
\node[circle, draw, inner sep=0pt, minimum size=7pt] at (t010) {};
\node[circle, draw, inner sep=0pt, minimum size=7pt] at (t100) {};
\node[circle, draw, inner sep=0pt, minimum size=7pt] at (t111) {};

\end{tikzpicture}
\caption{Strong tree envelopes for these pairs of rationals}
\end{figure}
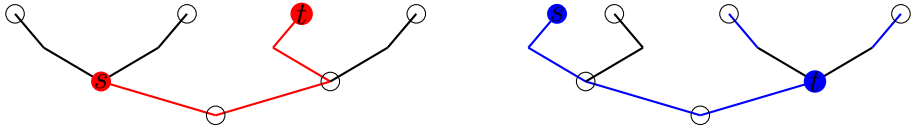

Trees can be used in various ways to code   structures.
As mentioned above, the nodes in the tree  $2^{<\om}$  ordered by   $\triangleleft$ produces a copy of the rationals.
In this setting, Sierpi\'{n}ski's coloring can be seen structurally as giving color red  to  pairs of  nodes  $s,t\in 2^{<\om}$ with $|s|\le |t|$ and $s\, \triangleleft\, t$,
and color blue otherwise.
Figure 1.\ gives pairs of nodes which are colored
red  (left) and blue (right).
These two configurations are examples of  {\em strong  similarity types},
 and  both of them persist in any subcopy of the rationals.
Figure 2.\ gives examples of strong trees which {\em envelope} such pairs of nodes.

The proof that  $T(2,\mathbb{Q})\le 2$ goes roughly as follows:
Consider pairs of nodes $s,t$ of different lengths  such that
 $s(|t|)=0$ if $|s|>|t|$, and $t(|s|)=0$ if $|t|>|s|$.
The two {\em strong similarity types} for such pairs
are seen in Figure 1.
Given a coloring of  pairs of nodes in $2^{<\om}$ into finitely many colors,
fix the strong similarity type on the left in Figure 1.
For each such pair $\{s,t\}$ in $2^{<\om}$ with this strong similarity type,
 there are finitely many $3$-strong trees enveloping  $s$ and $t$.
Give these  $3$-strong trees the coloring of the pair  which they envelope.
Apply Milliken's theorem to this coloring of $3$-strong trees to obtain an infinite  strong subtree $T$  where all pairs with this   strong similarity type have the same color.
Then do the same for the second strong similarity type,
obtaining an infinite strong subtree $T'$ of $T$.
Since the rationals are coded in any infinite strong subtree of $2^{<\om}$, one can now pull out an antichain of nodes, each of which has passing number $0$ at the levels of the other nodes, so that under the ordering $\triangleleft$, the nodes in this antichain represent $\mathbb{Q}$.
This produces a copy of the rationals in which all pairsets have at most two colors.


A similar but more involved strategy is behind the  finite big Ramsey degrees of the Rado graph.
For graphs,  by interpreting the lexicographic order on ordered graphs  due to  \Erdos, Hajnal, and P\'{o}sa in \cite{Erdos/Hajnal/Posa75} in terms of trees,
Sauer  showed in \cite{Sauer06}  that  nodes in the binary tree can code graphs.
Let $\mathbf{A}$ be a graph.
Enumerate the vertices of $\mathbf{A}$  as $\lgl v_n:n<N\rgl$.
A set of nodes $\{t_n:n<N\}$ in $2^{<\om}$ codes $\mathbf{A}$ if and only if
for each pair $m<n<N$,
\begin{equation}
v_n\ E\ v_m \Leftrightarrow  t_n(|t_m|)=1.
\end{equation}
The number $t_n(|t_m|)$ is called the {\em passing number} of $t_n$ at $t_m$.
See Figure 3.\ where the nodes $t_0,t_1,t_2$ in the tree code the graph $v_0,v_1,v_2$, which is a path of length two.

The Rado graph, $\mathcal{R}=(R,E)$, is the \Fraisse\ limit of the class of all finite graphs; as such, it is ultrahomogeneous and universal for all finite graphs.
\Erdos, Hajnal, and \Posa\ launched the investigation of finite big Ramsey degrees of the Rado graph in 1975,  when they showed that
there is a  2-coloring of edges such that every subcopy of the Rado graph retains both colors
 \cite{Erdos/Hajnal/Posa75},
 reminiscent of Sierpi\'{n}ski's result for pairs of rationals.
 Two decades later, Pouzet and Sauer proved that  for
each coloring of edges of $\mathcal{R}$ into finitely many colors, there is a subgraph $\mathcal{R}'$ isomorphic to $\mathcal{R}$ in which the edges have at most two colors
 \cite{Pouzet/Sauer96}.
In 2006, papers of
Sauer  \cite{Sauer06}  and of  Laflamme, Sauer, and Vuksanovic  \cite{Laflamme/Sauer/Vuksanovic06} combined to the exact big Ramsey degrees for all finite graphs within the Rado graph.
 In fact, these papers find big Ramsey degrees for a collection of binary relational structures, including the random tournament.
 The strategy is  similar to that outlined above for the rationals, but now the passing numbers at nodes code the edge/non-edge relation.


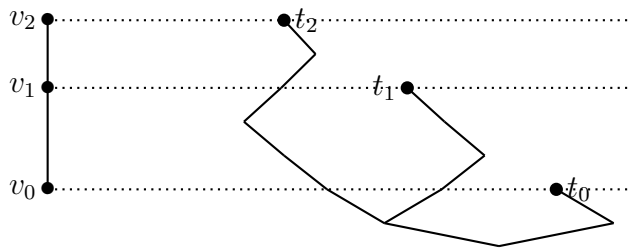
\begin{figure}\label{fig.bT}
\begin{tikzpicture}[grow'=up,scale=.3]
\tikzstyle{level 1}=[sibling distance=4in]
\tikzstyle{level 2}=[sibling distance=2in]
\tikzstyle{level 3}=[sibling distance=1.5in]
\tikzstyle{level 4}=[sibling distance=1.4in]
\tikzstyle{level 5}=[sibling distance=1.3in]
\tikzstyle{level 6}=[sibling distance=1.2in]
\tikzstyle{level 7}=[sibling distance=1.1in]
\tikzstyle{level 8}=[sibling distance=1in]
\node  {}
child{coordinate (t0)edge from parent[thick]
			child{coordinate (t00)
child{coordinate (t000)
child{coordinate (t0000)edge from parent[draw=black]
child{coordinate (t00000) edge from parent[draw=none]
}
child{coordinate (t00001)
child{edge from parent[draw=none]  coordinate (t000010)}
child{coordinate (t000011)
child{coordinate (t0000110)  edge from parent[draw=black]
}
child{edge from parent[draw=none]  coordinate (t0000111)}
}
}
}
child{edge from parent[draw=none] coordinate (t0001)}
}
child{ edge from parent[draw=none]coordinate(t001)}}
			child{ coordinate(t01)
child{ edge from parent[draw=none] coordinate(t010)}
child{coordinate(t011)
child{coordinate(t0110) edge from parent[draw=black]
child{coordinate(t01100)
}
child{ edge from parent[draw=none]  coordinate(t01101)}
}
child{ edge from parent[draw=none]  coordinate(t0111)}
}}}
		child{ coordinate(t1)edge from parent[thick]
			child{ coordinate(t10)
}
child{edge from parent[draw=none]   coordinate(t11)
} };

\node[right] at (t10) {$t_0$};
\node[left] at (t01100) {$t_1$};
\node[right] at (t0000110) {$t_2$};

\node[circle, fill=black,inner sep=0pt, minimum size=5pt] at (t10) {};
\node[circle, fill=black,inner sep=0pt, minimum size=5pt] at (t01100) {};
\node[circle, fill=black,inner sep=0pt, minimum size=5pt] at (t0000110) {};

\draw[thick, dotted] let \p1=(t10) in (-20,\y1) node (v0) {$\bullet$} -- (6,\y1);
\draw[thick, dotted] let \p1=(t01100) in (-20,\y1) node (v1) {$\bullet$} -- (6,\y1);
\draw[thick, dotted] let \p1= (t0000110) in (-20,\y1) node (v2) {$\bullet$} -- (6,\y1);

\node[left] at (v0) {$v_0$};
\node[left] at (v1) {$v_1$};
\node[left] at (v2) {$v_2$};

\draw[thick] (v0.center) to (v1.center) to (v2.center);

\end{tikzpicture}
\caption{A path of length two  coded by the nodes $t_0,t_1,t_2$}
\end{figure}

In \cite{Kechris/Pestov/Todorcevic05},
Kechris, Pestov, and Todorcevic asked,
Which structures have finite big Ramsey degrees?
In tandem with the results in \cite{Sauer06} and   \cite{Laflamme/Sauer/Vuksanovic06}, this
  sparked a new wave of  interest   in  big Ramsey degrees of \Fraisse\ structures.
As part of  his PhD work, Nguyen Van Th\'{e} proved  that
the countable
  ultrametric Urysohn space  with
any finite distance set
    has finite big Ramsey degrees \cite{NVT08}.
    While this result used Ramsey's theorem,
    the next result required a new extended version of Milliken's Theorem.
 Laflamme, Nguyen Van Th\'{e}, and Sauer
proved in \cite{Laflamme/NVT/Sauer10} that the structures
$\mathbb{Q}_n$ have finite big Ramsey degrees, where $\mathbb{Q}_n$ is
 the rationals as a linear order with an equivalence relation with $n$ equivalence classes, each of which is dense in $\mathbb{Q}$.
 See \cite{NVTHabil} for an excellent exposition of these and related results.

As for the Henson graphs,
Sauer proved in 1998 that the big Ramsey degree for edges in the triangle-free Henson graph is two  \cite{Sauer98}.
Further results were slow in coming, mainly because of lack of analogues for $k$-clique free graphs of the following two fortunate facts related to the Rado graph.
First,
the graph $\mathcal{B}$ induced by the countable  binary tree $2^{<\om}$ is universal for countable graphs.
Precisely, let
the vertices of  $\mathcal{B}$ be  the nodes in the tree $2^{<\om}$.
Define two nodes  to  have an edge between them  in $\mathcal{B}$ precisely when one of the nodes is longer than the other  and its passing number at the level of the shorter node is one.
Second,
Milliken's theorem for strong subtrees of $2^{<\om}$ can be applied to  subcopies of $\mathcal{B}$, and each infinite strong subtree  again  codes a graph which is universal for countable graphs.
After finitely many applications of
 Milliken's theorem to strong tree envelopes of finite antichains coding copies of  a given finite graph, one can  take a copy of the Rado graph where the big Ramsey degree of the finite graph under investigation is  finite.

 As analogues of these two facts were unknown for Henson graphs,
 the  two main themes in the work of \cite{DobrinenJML20} and \cite{DobrinenH_k19} were, first, to find means for coding the Henson graphs into trees in a way that the trees behaved  like strong trees, and
second, to prove  Ramsey theory for such trees.
In the next section, we discuss the
  methods  developed in these two papers to handle forbidden $k$-cliques.
These methods seem to be robust enough to handle many types of structures, including, perhaps surprisingly, infinite dimensional Ramsey theory of the Rado graph  (see \cite{DobrinenRado19}), which will also be discussed in the next section.

We close this section by pointing out the recent work of
\Masulovic\  in \cite{Masulovic18} which uses category theory to develop transfer principals for big Ramsey degrees.
In that paper, \Masulovic\
 widened the  investigation   of  big Ramsey degrees to    universal structures,  regardless of ultrahomogeneity.
 He has transferred several of the aforementioned results to prove finite big Ramsey degrees for some  new structures;
most strikingly,  these include \Fraisse\ limits of classes of finite metric spaces with finite distance sets satisfying certain properties.


\section{Trees with coding nodes, their Ramsey theory,\\ and applications to Ramsey theory of infinite graphs}\label{sec.sct}

In order to  investigate  big Ramsey degrees of the Henson graphs,
the first task was to find some  way of representing
 $k$-clique-free graphs via trees.
 Since $k$-cliques are forbidden, there needs to be
 some method for determining
which nodes should be allowed to split
and which ones should not, so that
$\mathcal{H}_k$  is  represented   while the  splitting  is `maximal' in some sense.
In particular, the trees need to be perfect in order for Ramsey theory to have any chance of development.

To achieve this,
in the construction of such trees, certain nodes are distinguished to code certain vertices of $\mathcal{H}_k$.
 This way, one keeps track of the finite graph that is already coded, ensuring that one knows how to branch maximally, subject to  never coding a $k$-clique in the future.
This led to  the following notions of trees with coding nodes and  strong coding trees.


\subsection{Strong coding trees}\label{subsec.sct}

The following definitions and theorems  are taken from \cite{DobrinenJML20} and \cite{DobrinenH_k19}.
A {\em tree with coding nodes} is a structure $\lgl T,N; \sse, <,c^T\rgl$  in the language $\mathcal{L}=\{\sse,<,c\}$,
where $\sse$ and $<$ are binary relation symbols and $c$ is a unary function symbol satisfying the following:
$T\sse 2^{<\om}$ and $(T,\sse)$ is a tree.
$N\le \om$ and $<$ is the standard linear order on $N$.
$c^T:N\ra T$ is injective, and $m<n<N \longrightarrow |c^T(m)|<|c^T(n)|$.
$c^T(n)$ is the {\em $n$-th coding node in $T$},  and is usually denoted $c_n^T$.

Notice that a  collection of   coding  nodes
$\{c^T_{n_i}:i<k\}$ in T
 {\em  codes a $k$-clique}
  if and only if
  $c^T_{n_j}(|c^T_{n_i}|)=1$ for all  $i<j<k$.
To ensure that a tree never codes a $k$-clique,   the following branching criterion is introduced.
We say that a tree $T$ with coding nodes $\lgl c^T_n:n<N\rgl$ satisfies the {\em  $K_k$-Free Branching  Criterion ($k$-FBC)}
if  for each non-maximal node $t\in T$,
$t^{\frown}0$ is always  in $T$,
and
$t^{\frown}1$ is in $T$ if and only if  any coding node extending  $t^{\frown}1$  cannot code a $k$-clique with coding nodes in $T$ of  length  less than or equal to the length of $t$.
It is a useful fact  that given any tree $T$ with coding nodes and no maximal nodes
 satisfying the $k$-FBC, and  in which the set of  coding nodes are dense,
the coding nodes in $T$  code $\mathcal{H}_k$
(Theorem 4.9, \cite{DobrinenH_k19}).


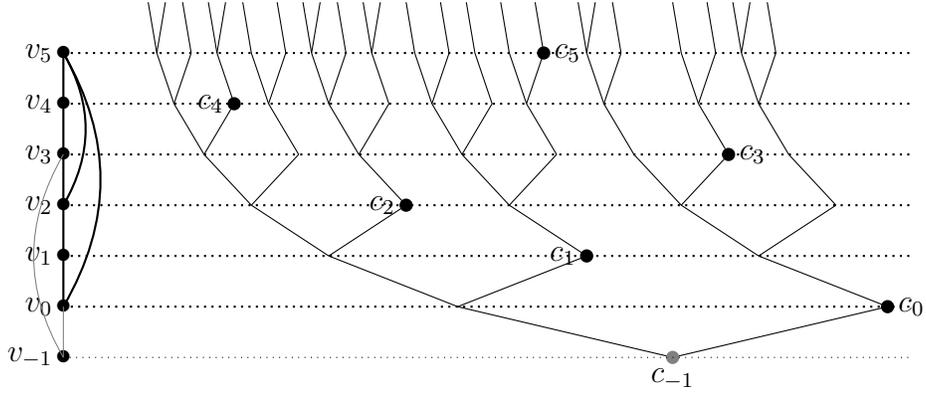
\begin{figure}\label{fig.bS3}
\begin{tikzpicture}[grow'=up,scale=.45]

\tikzstyle{level 1}=[sibling distance=5in]
\tikzstyle{level 2}=[sibling distance=3in]
\tikzstyle{level 3}=[sibling distance=1.8in]
\tikzstyle{level 4}=[sibling distance=1.1in]
\tikzstyle{level 5}=[sibling distance=.7in]
\tikzstyle{level 6}=[sibling distance=0.4in]
\tikzstyle{level 7}=[sibling distance=0.2in]

\node {} coordinate(t)
child{coordinate (t0)
			child{coordinate (t00)
child{coordinate (t000)
child {coordinate(t0000)
child{coordinate(t00000)
child{coordinate(t000000)
child{coordinate(t0000000)}
child{coordinate(t0000001)}
}
child{coordinate(t000001)
child{coordinate(t0000010)}
child{ edge from parent[draw=none]  coordinate(t0000011)}
}
}
child{coordinate(t00001)
child{coordinate(t000010)
child{coordinate(t0000100)}
child{coordinate(t0000101)}
}
child{ edge from parent[draw=none]  coordinate(t000011)
}
}}
child {coordinate(t0001)
child {coordinate(t00010)
child {coordinate(t000100)
child {coordinate(t0001000)}
child { edge from parent[draw=none] coordinate(t0001001)}
}
child {coordinate(t000101)
child {coordinate(t0001010)}
child { edge from parent[draw=none]  coordinate(t0001011)}
}
}
child{coordinate(t00011) edge from parent[draw=none] }}}
child{ coordinate(t001)
child{ coordinate(t0010)
child{ coordinate(t00100)
child{ coordinate(t001000)
child{ coordinate(t0010000)}
child{ coordinate(t0010001)}
}
child{ coordinate(t001001)
child{ coordinate(t0010010)}
child{ edge from parent[draw=none] coordinate(t0010011)}
}
}
child{ coordinate(t00101)
child{ coordinate(t001010)
child{ coordinate(t0010100)}
child{ coordinate(t0010101)}
}
child{   edge from parent[draw=none]coordinate(t001011)
}
}}
child{  edge from parent[draw=none] coordinate(t0011)}}}
			child{ coordinate(t01)
child{ coordinate(t010)
child{ coordinate(t0100)
child{ coordinate(t01000)
child{ coordinate(t010000)
child{ coordinate(t0100000)}
child{ edge from parent[draw=none]  coordinate(t0100001)}
}
child{ coordinate(t010001)
child{ coordinate(t0100010)}
child{edge from parent[draw=none]  coordinate(t0100011)}
}
}
child{ coordinate(t01001)
child{ coordinate(t010010)
child{ coordinate(t0100100)}
child{ edge from parent[draw=none]  coordinate(t0100101)}
}
child{edge from parent[draw=none]  coordinate(t010011)}
}}
child{ coordinate(t0101)
child{ coordinate(t01010)
child{ coordinate(t010100)
child{ coordinate(t0101000)}
child{edge from parent[draw=none]  coordinate(t0101001)}
}
child{ coordinate(t010101)
child{ coordinate(t0101010)}
child{edge from parent[draw=none]  coordinate(t0101011)}
}
}
child{  edge from parent[draw=none]  coordinate(t01011)
}}}
child{ edge from parent[draw=none]  coordinate(t011)}}}
		child{ coordinate(t1)
			child{ coordinate(t10)
child{ coordinate(t100)
child{ coordinate(t1000)
child{ coordinate(t10000)
child{ coordinate(t100000)
child{ coordinate(t1000000)}
child{ coordinate(t1000001)}
}
child{ coordinate(t100001)
child{ coordinate(t1000010)}
child{ edge from parent[draw=none] coordinate(t1000011)}
}
}
child{ edge from parent[draw=none] coordinate(t10001)
}}
child{ coordinate(t1001)
child{ coordinate(t10010)
child{ coordinate(t100100)
child{ coordinate(t1001000)}
child{edge from parent[draw=none]  coordinate(t1001001)}
}
child{ coordinate(t100101)
child{ coordinate(t1001010)}
child{edge from parent[draw=none]   coordinate(t1001011)}
}
}
child{  edge from parent[draw=none] coordinate(t10011)
}}}
child{ coordinate(t101)
child{ coordinate(t1010)
child{ coordinate(t10100)
child{ coordinate(t101000)
child{ coordinate(t1010000)}
child{ coordinate(t1010001)}
}
child{ coordinate(t101001)
child{ coordinate(t1010010)}
child{   edge from parent[draw=none]   coordinate(t1010011)}
}
}
child{   edge from parent[draw=none]  coordinate(t10101)
}}
child{ edge from parent[draw=none] coordinate(t1011)}}}
child{  edge from parent[draw=none] coordinate(t11)} };

\node[below] at (t) {$c_{-1}$};
\node[right] at (t1) {$c_0$};
\node[left] at (t01) {$c_1$};
\node[left] at (t001) {$c_2$};
\node[right] at (t1001) {$c_3$};
\node[left] at (t00001) {$c_4$};
\node[right] at (t010101) {$c_5$};

\node[circle, fill=gray,inner sep=0pt, minimum size=5pt] at (t) {};
\node[circle, fill=black,inner sep=0pt, minimum size=5pt] at (t1) {};
\node[circle, fill=black,inner sep=0pt, minimum size=5pt] at (t01) {};
\node[circle, fill=black,inner sep=0pt, minimum size=5pt] at (t001) {};
\node[circle, fill=black,inner sep=0pt, minimum size=5pt] at (t1001) {};
\node[circle, fill=black,inner sep=0pt, minimum size=5pt] at (t00001) {};
\node[circle, fill=black,inner sep=0pt, minimum size=5pt] at (t010101) {};

\draw[dotted] let \p1=(t) in (-18,\y1) node (v00) {$\bullet$} -- (7,\y1);
\draw[thick, dotted] let \p1=(t1) in (-18,\y1) node (v0) {$\bullet$} -- (7,\y1);
\draw[thick, dotted] let \p1=(t01) in (-18,\y1) node (v1) {$\bullet$} -- (7,\y1);
\draw[thick, dotted] let \p1=(t001) in (-18,\y1) node (v2) {$\bullet$} -- (7,\y1);
\draw[thick, dotted] let \p1=(t1001) in (-18,\y1) node (v3) {$\bullet$} -- (7,\y1);
\draw[thick, dotted] let \p1=(t00001) in (-18,\y1) node (v4) {$\bullet$} -- (7,\y1);
\draw[thick, dotted] let \p1=(t010101) in (-18,\y1) node (v5) {$\bullet$} -- (7,\y1);

\node[left] at (v00) {$v_{-1}$};
\node[left] at (v0) {$v_0$};
\node[left] at (v1) {$v_1$};
\node[left] at (v2) {$v_2$};
\node[left] at (v3) {$v_3$};
\node[left] at (v4) {$v_4$};
\node[left] at (v5) {$v_5$};

\draw[thick] (v0.center) to (v1.center) to (v2.center) to (v3.center);
\draw[thick] (v3.center) to (v4.center) to (v5.center);
\draw[thick] (v0.center) to [bend right] (v5.center);
\draw[thick] (v5.center) to [bend left] (v2.center);
\draw[gray] (v00.center) to [bend left] (v3.center);
\draw[gray] (v00.center) to (v0.center);

\end{tikzpicture}
\caption{A strong triangle-free tree $\bS_3$ densely coding $\mathcal{H}_3$}
\end{figure}


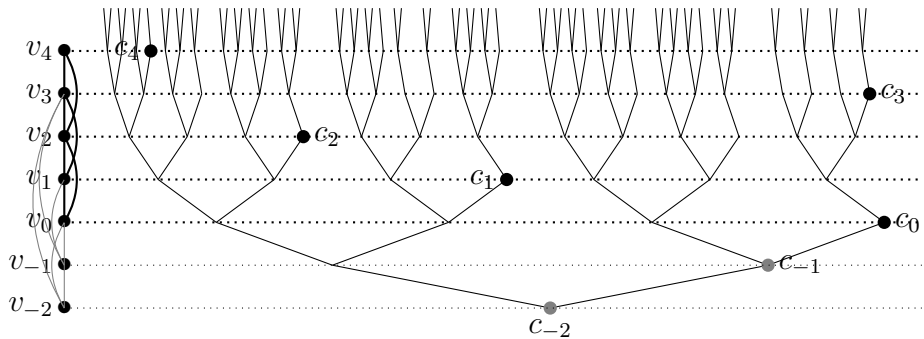
\begin{figure}
\begin{tikzpicture}[grow'=up,scale=.38]

\tikzstyle{level 1}=[sibling distance=6in]
\tikzstyle{level 2}=[sibling distance=3.2in]
\tikzstyle{level 3}=[sibling distance=1.6in]
\tikzstyle{level 4}=[sibling distance=.8in]
\tikzstyle{level 5}=[sibling distance=.4in]
\tikzstyle{level 6}=[sibling distance=0.2in]
\tikzstyle{level 7}=[sibling distance=0.1in]

\node {} coordinate(t)
child{coordinate (t0)
			child{coordinate (t00)
child{coordinate (t000)
child {coordinate(t0000)
child{coordinate(t00000)
child{coordinate(t000000)
child{coordinate(t0000000)}
child{coordinate(t0000001)}
}
child{coordinate(t000001)
child{coordinate(t0000010)}
child{ coordinate(t0000011)}
}
}
child{coordinate(t00001)
child{coordinate(t000010)
child{coordinate(t0000100)}
child{coordinate(t0000101)}
}
child{ coordinate(t000011)
child{ coordinate(t0000110)}
child{ coordinate(t0000111) edge from parent[draw=none] }
}
}}
child {coordinate(t0001)
child {coordinate(t00010)
child {coordinate(t000100)
child {coordinate(t0001000)}
child {  coordinate(t0001001)}
}
child {coordinate(t000101)
child {coordinate(t0001010)}
child { coordinate(t0001011)}
}
}
child{coordinate(t00011)
child{coordinate(t000110)
child{coordinate(t0001100)}
child{coordinate(t0001101)}
 }
child{coordinate(t000111)edge from parent[draw=none] }
 }}}
child{ coordinate(t001)
child{ coordinate(t0010)
child{ coordinate(t00100)
child{ coordinate(t001000)
child{ coordinate(t0010000)}
child{ coordinate(t0010001)}
}
child{ coordinate(t001001)
child{ coordinate(t0010010)}
child{ coordinate(t0010011)}
}
}
child{ coordinate(t00101)
child{ coordinate(t001010)
child{ coordinate(t0010100)}
child{ coordinate(t0010101)}
}
child{ coordinate(t001011)
child{ coordinate(t0010110)}
child{ coordinate(t0010111) edge from parent[draw=none] }
}
}}
child{  coordinate(t0011)
child{  coordinate(t00110)
child{  coordinate(t001100)
child{  coordinate(t0011000)}
child{  coordinate(t0011001)}
}
child{  coordinate(t001101)
child{  coordinate(t0011010)}
child{  coordinate(t0011011)}
}
}
child{  coordinate(t00111)edge from parent[draw=none]}
}}}
			child{ coordinate(t01)
child{ coordinate(t010)
child{ coordinate(t0100)
child{ coordinate(t01000)
child{ coordinate(t010000)
child{ coordinate(t0100000)}
child{  coordinate(t0100001)}
}
child{ coordinate(t010001)
child{ coordinate(t0100010)}
child{ coordinate(t0100011)}
}
}
child{ coordinate(t01001)
child{ coordinate(t010010)
child{ coordinate(t0100100)}
child{ coordinate(t0100101)}
}
child{ coordinate(t010011)
child{ coordinate(t0100110)}
child{ coordinate(t0100111)edge from parent[draw=none]  }
}
}}
child{ coordinate(t0101)
child{ coordinate(t01010)
child{ coordinate(t010100)
child{ coordinate(t0101000)}
child{ coordinate(t0101001)}
}
child{ coordinate(t010101)edge from parent[draw=none]}
}
child{  coordinate(t01011)
child{  coordinate(t010110)
child{  coordinate(t0101100)}
child{  coordinate(t0101101)}
}
child{  coordinate(t010111)edge from parent[draw=none]  }
}}}
child{ coordinate(t011)
child{ coordinate(t0110)
child{ coordinate(t01100)
child{ coordinate(t011000)
child{ coordinate(t0110000)}
child{ coordinate(t0110001)}
}
child{ coordinate(t011001)
child{ coordinate(t0110010)}
child{ coordinate(t0110011)}
}
}
child{ coordinate(t01101)
child{ coordinate(t011010)
child{ coordinate(t0110100)}
child{ coordinate(t0110101)}
}
child{ coordinate(t011011)
child{ coordinate(t0110110)}
child{ coordinate(t0110111)edge from parent[draw=none] }}
}
}
child{ coordinate(t0111)edge from parent[draw=none]
}
}}}
		child{ coordinate(t1)
			child{ coordinate(t10)
child{ coordinate(t100)
child{ coordinate(t1000)
child{ coordinate(t10000)
child{ coordinate(t100000)
child{ coordinate(t1000000)}
child{ coordinate(t1000001)}
}
child{ coordinate(t100001)
child{ coordinate(t1000010)}
child{ coordinate(t1000011)}
}
}
child{ coordinate(t10001)
child{ coordinate(t100010)
child{ coordinate(t1000100)}
child{ coordinate(t1000101)}
}
child{ coordinate(t100011)
child{ coordinate(t1000110)}
child{ coordinate(t1000111)edge from parent[draw=none]  }
}
}}
child{ coordinate(t1001)
child{ coordinate(t10010)
child{ coordinate(t100100)
child{ coordinate(t1001000)}
child{ coordinate(t1001001)}
}
child{ coordinate(t100101)
child{ coordinate(t1001010)}
child{ coordinate(t1001011)}
}
}
child{  coordinate(t10011)
child{  coordinate(t100110)
child{  coordinate(t1001100)}
child{  coordinate(t1001101)}
}
child{  coordinate(t100111)edge from parent[draw=none] }
}}}
child{ coordinate(t101)
child{ coordinate(t1010)
child{ coordinate(t10100)
child{ coordinate(t101000)
child{ coordinate(t1010000)}
child{ coordinate(t1010001)}
}
child{ coordinate(t101001)
child{ coordinate(t1010010)}
child{   coordinate(t1010011)}
}
}
child{  coordinate(t10101)
child{  coordinate(t101010)
child{  coordinate(t1010100)}
child{  coordinate(t1010101)}
}
child{  coordinate(t101011)
child{  coordinate(t1010110)}
child{  coordinate(t1010111)edge from parent[draw=none]}
}
}}
child{ coordinate(t1011)
child{ coordinate(t10110)
child{ coordinate(t101100)
child{ coordinate(t1011000)}
child{ coordinate(t1011001)}
}
child{ coordinate(t101101)
child{ coordinate(t1011010)}
child{ coordinate(t1011011)}
}
}
child{ coordinate(t10111)
edge from parent[draw=none]  }
}}}
child{ coordinate(t11)
child{coordinate (t110)
child{coordinate (t1100)
child{coordinate (t11000)
child{coordinate(t110000)
child{coordinate(t1100000)}
child{coordinate(t1100001)}
}
child{coordinate(t110001) edge from parent[draw=none]  }}
child{coordinate(t11001)
child{coordinate(t110010)
child{coordinate(t1100100)}
child{coordinate(t1100101)}}
child{coordinate(t110011) edge from parent[draw=none]  }
}}
child{coordinate (t1101)
child{coordinate(t11010)
child{coordinate(t110100)
child{coordinate(t1101000)}
child{coordinate(t1101001)}}
child{coordinate(t110101) edge from parent[draw=none]  }
}
child{coordinate(t11011)
child{coordinate(t110110)
child{coordinate(t1101100)}
child{coordinate(t1101101)}
}
child{coordinate(t110111) edge from parent[draw=none]  }
}}}
child{coordinate(t111) edge from parent[draw=none] }} };

\node[below] at (t) {$c_{-2}$};
\node[right] at (t1) {$c_{-1}$};
\node[right] at (t11) {$c_0$};
\node[left] at (t011) {$c_1$};
\node[right] at (t0011) {$c_2$};
\node[right] at (t11011) {$c_3$};
\node[left] at (t000011) {$c_4$};

\node[circle, fill=gray,inner sep=0pt, minimum size=5pt] at (t) {};
\node[circle, fill=gray,inner sep=0pt, minimum size=5pt] at (t1) {};
\node[circle, fill=black,inner sep=0pt, minimum size=5pt] at (t11) {};
\node[circle, fill=black,inner sep=0pt, minimum size=5pt] at (t011) {};
\node[circle, fill=black,inner sep=0pt, minimum size=5pt] at (t0011) {};
\node[circle, fill=black,inner sep=0pt, minimum size=5pt] at (t11011) {};
\node[circle, fill=black,inner sep=0pt, minimum size=5pt] at (t000011) {};

\draw[dotted] let \p1=(t) in (-17,\y1) node (v02) {$\bullet$} -- (13,\y1);
\draw[dotted] let \p1=(t1) in (-17,\y1) node (v01) {$\bullet$} -- (13,\y1);
\draw[thick, dotted] let \p1=(t11) in (-17,\y1) node (v0) {$\bullet$} -- (13,\y1);
\draw[thick, dotted] let \p1=(t011) in (-17,\y1) node (v1) {$\bullet$} -- (13,\y1);
\draw[thick, dotted] let \p1=(t0011) in (-17,\y1) node (v2) {$\bullet$} -- (13,\y1);
\draw[thick, dotted] let \p1=(t11011) in (-17,\y1) node (v3) {$\bullet$} -- (13,\y1);
\draw[thick, dotted] let \p1=(t000011) in (-17,\y1) node (v4) {$\bullet$} -- (13,\y1);

\node[left] at (v02) {$v_{-2}$};
\node[left] at (v01) {$v_{-1}$};
\node[left] at (v0) {$v_0$};
\node[left] at (v1) {$v_1$};
\node[left] at (v2) {$v_2$};
\node[left] at (v3) {$v_3$};
\node[left] at (v4) {$v_4$};

\draw[gray] (v02.center) to (v01.center) to (v0.center) to [bend right] (v02.center);
\draw[thick] (v0.center) to (v1.center) ;
\draw[gray] (v1.center) to [bend right] (v01.center);
\draw[thick] (v1.center) to (v2.center)to [bend left] (v0.center) ;
\draw[thick] (v2.center) to (v3.center)to [bend left] (v1.center) ;
\draw[gray] (v3.center) to [bend right] (v01.center);
\draw[gray] (v3.center) to [bend right] (v02.center);
\draw[thick] (v3.center) to (v4.center)to [bend left] (v2.center) ;

\end{tikzpicture}
\caption{A strong $K_4$-free tree $\bS_4$ densely coding $\mathcal{H}_4$}
\end{figure}


The trees $\bS_3$ and $\bS_4$ in
Figures 4.\ and 5.\ have coding nodes which
 code Henson graphs $\mathcal{H}_3$ and $\mathcal{H}_4$, respectively.
These trees feature  the main structural ideas behind strong coding trees;
one can extrapolate to envisage  $\bS_k$, for any $k\ge 5$ (precise constructions are given in \cite{DobrinenH_k19}).
The gray nodes $c_{-(k-2)},\dots,c_{-1}$
are pseudo-coding nodes, which code a $(k-1)$-clique.
They help  to set up the tree structure so that subtrees  of $\mathbb{S}_k$  coding $\mathcal{H}_k$ can be isomorphic to $\mathbb{S}_k$.
The vertex $v_0$ is to be thought of as
 forming a $(k-1)$-clique with some  vertices in a larger ambient copy of $\mathcal{H}_k$.
This has the effect that each coding node in $\bS_k$
does not split.

Perhaps the most illustrative way of thinking of these trees is the following:
The nodes  in the tree the level of $c_n$ in the tree are coding all finite partial types over
the graph on vertices $\{v_{-(k-2)},\dots, v_{n-1}\}$.
In this way, a strong coding tree is really just  a means for visualizing the  finite partial types over an  (ordered) initial segment of the graph $\mathcal{H}_k$.

The one small but insurmountable catch to these trees is that  having coding nodes and splitting nodes at the same levels prevents  the development of Ramsey theory on subtrees isomorphic to $\mathbb{S}_k$.
Ironically, the failure occurs at the most basic  level:
 There are bad colorings of coding nodes for which no subtree isomorphic to $\mathbb{S}_k$ has one color (see Example 3.18 in \cite{DobrinenJML20}).
There turns out to be a simple solution:
Skew the trees so that each level of the tree has at most one coding node or splitting node, but never both.
Let   $\mathbb{T}_k$ denote the skewed version of $\bS_k$ (see Figure 6.\ for $\mathbb{T}_3$).


\begin{figure}\label{fig.bT}
\begin{tikzpicture}[grow'=up,scale=.3]
\tikzstyle{level 1}=[sibling distance=7in]
\tikzstyle{level 2}=[sibling distance=2in]
\tikzstyle{level 3}=[sibling distance=2in]
\tikzstyle{level 4}=[sibling distance=2in]
\tikzstyle{level 5}=[sibling distance=1in]
\tikzstyle{level 6}=[sibling distance=1.8in]
\tikzstyle{level 7}=[sibling distance=.8in]
\tikzstyle{level 8}=[sibling distance=.8in]
\tikzstyle{level 9}=[sibling distance=.8in]
\tikzstyle{level 10}=[sibling distance=1.2in]
\tikzstyle{level 11}=[sibling distance=.3in]
\tikzstyle{level 12}=[sibling distance=.3in]
\tikzstyle{level 13}=[sibling distance=.3in]
\tikzstyle{level 14}=[sibling distance=.5in]

\node {} coordinate(t)

child{coordinate (t0)
			child{coordinate (t00)
child{coordinate (t000)
child{coordinate (t0000)
child{coordinate (t00000)
child{coordinate (t000000)
child{coordinate (t0000000)
child{coordinate (t00000000)
child{coordinate (t000000000)
child{coordinate (t0000000000)
child{coordinate (t00000000000)
child{coordinate (t000000000000)
child{coordinate (t0000000000000)
child{coordinate (t00000000000000)}
child{edge from parent[draw=none] coordinate (t00000000000001)}
}
child{coordinate (t0000000000001)
child{edge from parent[draw=none]  coordinate (t00000000000010)}
child{coordinate (t00000000000011)}
}
}
child{edge from parent[draw=none] coordinate (t000000000001)}
}
child{edge from parent[draw=none]  coordinate (t00000000001)}
}
child{edge from parent[draw=none]  coordinate (t0000000001)}
}
child{coordinate (t000000001)
child{edge from parent[draw=none]  coordinate (t0000000010)}
child{coordinate (t0000000011)
child{coordinate (t00000000110)
child{coordinate (t000000001100)
child{coordinate (t0000000011000)
child{coordinate (t00000000110000)}
child{edge from parent[draw=none] coordinate (t00000000110001)}
}
child{edge from parent[draw=none] coordinate (t0000000011001)}
}
child{edge from parent[draw=none]  coordinate (t000000001101)}
}
child{edge from parent[draw=none]   coordinate (t00000000111)}
}
}
}
child{edge from parent[draw=none] coordinate (t00000001)}
}
child{edge from parent[draw=none] coordinate (t0000001)}
}
child{edge from parent[draw=none]  coordinate (t000001)}
}
child{coordinate (t00001)
child{edge from parent[draw=none]  coordinate (t000010)}
child{coordinate (t000011)
child{coordinate (t0000110)
child{coordinate (t00001100)
child{coordinate (t000011000)
child{coordinate (t0000110000)
child{coordinate (t00001100000)
child{coordinate (t000011000000)
child{coordinate (t0000110000000)
child{coordinate (t00001100000000)}
child{edge from parent[draw=none] coordinate (t00001100000001)}
}
child{edge from parent[draw=none] coordinate (t0000110000001)}
}
child{coordinate (t000011000001)
child{coordinate (t0000110000010)
child{edge from parent[draw=none]  coordinate (t00001100000100)}
child{coordinate (t00001100000101)}
}
child{edge from parent[draw=none] coordinate (t0000110000011)}
}
}
child{edge from parent[draw=none] coordinate (t00001100001)}
}
child{edge from parent[draw=none]  coordinate (t0000110001)}
}
child{edge from parent[draw=none] coordinate (t000011001)}
}
child{edge from parent[draw=none]  coordinate (t00001101)}
}
child{edge from parent[draw=none]  coordinate (t0000111)}
}
}
}
child{edge from parent[draw=none] coordinate (t0001)}
}
child{ edge from parent[draw=none]coordinate(t001)}}
			child{ coordinate(t01)
child{ edge from parent[draw=none] coordinate(t010)}
child{coordinate(t011)
child{coordinate(t0110)
child{coordinate(t01100)
child{coordinate(t011000)
child{coordinate(t0110000)
child{coordinate(t01100000)
child{coordinate(t011000000)
child{coordinate(t0110000000)
child{coordinate(t01100000000)
child{coordinate(t011000000000)
child{coordinate(t0110000000000)
child{coordinate(t01100000000000)}
child{edge from parent[draw=none] coordinate(t01100000000001)}
}
child{ edge from parent[draw=none]coordinate(t0110000000001)}
}
child{ edge from parent[draw=none]  coordinate(t011000000001)}
}
child{coordinate(t01100000001)
child{coordinate(t011000000010)
child{coordinate(t0110000000100)
child{edge from parent[draw=none] coordinate(t01100000001000)}
child{coordinate(t01100000001001)}
}
child{edge from parent[draw=none] coordinate(t0110000000101)}
}
child{edge from parent[draw=none]   coordinate(t011000000011)}
}
}
child{ edge from parent[draw=none]  coordinate(t0110000001)}
}
child{ edge from parent[draw=none] coordinate(t011000001)}
}
child{coordinate(t01100001)
child{coordinate(t011000010)
child{edge from parent[draw=none] coordinate(t0110000100)}
child{coordinate(t0110000101)
child{coordinate(t01100001010)
child{coordinate(t011000010100)
child{coordinate(t0110000101000)
child{coordinate(t01100001010000)}
child{edge from parent[draw=none] coordinate(t01100001010001)}
}
child{edge from parent[draw=none]  coordinate(t0110000101001)}
}
child{edge from parent[draw=none]  coordinate(t011000010101)}
}
child{edge from parent[draw=none] coordinate(t01100001011)}
}
}
child{edge from parent[draw=none]  coordinate(t011000011)}
}
}
child{ edge from parent[draw=none] coordinate(t0110001)}
}
child{ edge from parent[draw=none]  coordinate(t011001)}
}
child{ edge from parent[draw=none]  coordinate(t01101)}
}
child{ edge from parent[draw=none]  coordinate(t0111)}
}}}
		child{ coordinate(t1)
			child{ coordinate(t10)
child{ coordinate(t100)
child{ coordinate(t1000)
child{ coordinate(t10000)
child{ coordinate(t100000)
child{ coordinate(t1000000)
child{ coordinate(t10000000)
child{ coordinate(t100000000)
child{ coordinate(t1000000000)
child{ coordinate(t10000000000)
child{ coordinate(t100000000000)
child{ coordinate(t1000000000000)
child{ coordinate(t10000000000000)}
child{ edge from parent[draw=none]   coordinate(t10000000000001)}
}
child{edge from parent[draw=none]   coordinate(t1000000000001)}
}
child{ edge from parent[draw=none]   coordinate(t100000000001)}
}
child{edge from parent[draw=none]  coordinate(t10000000001)}
}
child{ edge from parent[draw=none]  coordinate(t1000000001)}
}
child{edge from parent[draw=none]  coordinate(t100000001)}
}
child{   edge from parent[draw=none] coordinate(t10000001)}
}
child{ coordinate(t1000001)
child{ coordinate(t10000010)
child{ coordinate(t100000100)
child{ edge from parent[draw=none]  coordinate(t1000001000)}
child{ coordinate(t1000001001)
child{ coordinate(t10000010010)
child{ coordinate(t100000100100)
child{ coordinate(t1000001001000)
child{ coordinate(t10000010010000)}
child{ edge from parent[draw=none] coordinate(t10000010010001)}
}
child{edge from parent[draw=none]   coordinate(t1000001001001)}
}
child{edge from parent[draw=none]  coordinate(t100000100101)}
}
child{ edge from parent[draw=none] coordinate(t10000010011)}
}
}
child{  edge from parent[draw=none] coordinate(t100000101)}
}
child{   edge from parent[draw=none] coordinate(t10000011)}
}
}
child{  edge from parent[draw=none] coordinate(t100001)}
}
child{ edge from parent[draw=none]  coordinate(t10001)}
}
child{ coordinate(t1001)
child{ coordinate(t10010)
child{ edge from parent[draw=none]   coordinate(t100100)}
child{ coordinate(t100101)
child{ coordinate(t1001010)
child{ coordinate(t10010100)
child{ coordinate(t100101000)
child{ coordinate(t1001010000)
child{ coordinate(t10010100000)
child{ coordinate(t100101000000)
child{ coordinate(t1001010000000)
child{ coordinate(t10010100000000)}
child{edge from parent[draw=none] coordinate(t10010100000001)}
}
child{ edge from parent[draw=none]coordinate(t1001010000001)}
}
child{edge from parent[draw=none]  coordinate(t100101000001)}
}
child{edge from parent[draw=none]  coordinate(t10010100001)}
}
child{ edge from parent[draw=none]  coordinate(t1001010001)}
}
child{ edge from parent[draw=none] coordinate(t100101001)}
}
child{edge from parent[draw=none]   coordinate(t10010101)}
}
child{ edge from parent[draw=none]   coordinate(t1001011)}
}
}
child{  edge from parent[draw=none] coordinate(t10011)}
}
}
child{edge from parent[draw=none]  coordinate(t101)}}
child{edge from parent[draw=none]   coordinate(t11)
} };

\node[below] at (t) {$d_0=c_{-1}$};
\node[left] at (t0) {$d_1$};
\node[right] at (t10) {$c_0$};
\node[left] at (t10) {$d_2$};
\node[left] at (t100) {$d_3$};
\node[left] at (t0000) {$d_4$};
\node[left] at (t01100) {$d_5$};
\node[right] at (t01100) {$c_1$};
\node[left] at (t100000) {$d_6$};
\node[left] at (t0110000) {$d_7$};
\node[left] at (t00000000) {$d_8$};
\node[right] at (t000011000) {$c_2$};
\node[left] at (t000011000) {$d_9$};

\node[left] at (t0110000000) {$d_{10}$};
\node[left] at (t00001100000) {$d_{11}$};
\node[left] at (t000000000000) {$d_{12}$};
\node[left] at (t1000001001000) {$d_{13}$};
\node[right] at (t1000001001000)  {$c_3$};

\node[circle, fill=gray,inner sep=0pt, minimum size=5pt] at (t) {};
\node[circle, fill=black,inner sep=0pt, minimum size=5pt] at (t10) {};
\node[circle, fill=black,inner sep=0pt, minimum size=5pt] at (t01100) {};
\node[circle, fill=black,inner sep=0pt, minimum size=5pt] at (t000011000) {};
\node[circle, fill=black,inner sep=0pt, minimum size=5pt] at (t1000001001000) {};

\draw[thick, dotted] let \p1=(t) in (-28,\y1) node (v01) {$\bullet$} -- (9,\y1);
\draw[thick, dotted] let \p1=(t10) in (-28,\y1) node (v0) {$\bullet$} -- (9,\y1);
\draw[thick, dotted] let \p1=(t01100) in (-28,\y1) node (v1) {$\bullet$} -- (9,\y1);
\draw[thick, dotted] let \p1= (t000011000) in (-28,\y1) node (v2) {$\bullet$} -- (9,\y1);
\draw[thick, dotted] let \p1=  (t1000001001000) in (-28,\y1) node (v3) {$\bullet$} -- (9,\y1);

\node[left] at (v01) {$v_{-1}$};
\node[left] at (v0) {$v_0$};
\node[left] at (v1) {$v_1$};
\node[left] at (v2) {$v_2$};
\node[left] at (v3) {$v_3$};

\draw[gray] (v0.center) to (v01.center)  to [bend left] (v3.center);
\draw[thick] (v0.center) to (v1.center) to (v2.center) to (v3.center) ;

\end{tikzpicture}
\caption{Strong $\mathcal{H}_3$-coding tree $\bT_3$}
\end{figure}
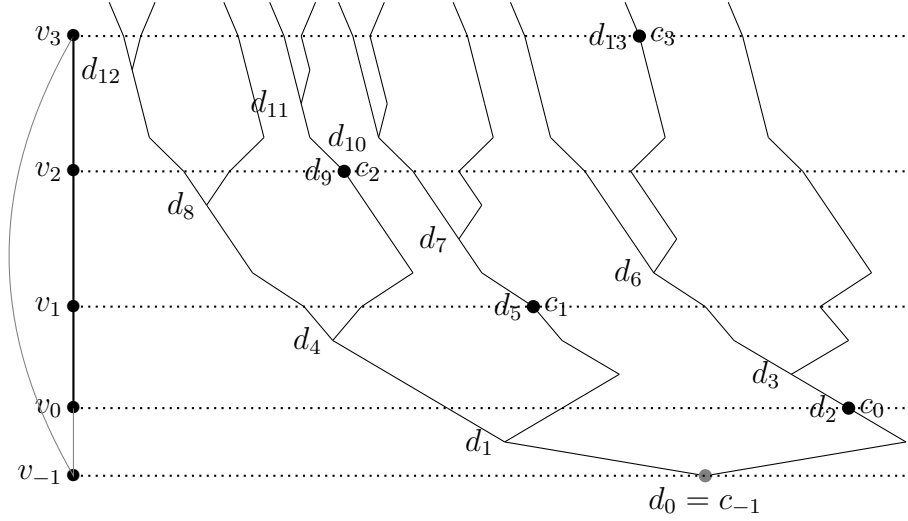

We  work with a collection $\mathcal{T}_k$ of infinite subtrees of
$\mathbb{T}_k$, each of which is isomorphic to $\mathbb{T}_k$ in a strong sense:
Let $k\ge 3$ be given and let  $S,T\sse \bT_k$ be meet-closed subsets.
A bijection $f:S\ra T$ is a {\em strong similarity map}  if for all nodes $s,t,u,v\in S$, the following hold:
\begin{enumerate}
\item
$f$ preserves lexicographic order.
\item
$f$ preserves meets, and hence splitting nodes.
\item
$f$ preserves relative lengths.
\item
$f$ preserves initial segments.
\item
$f$ preserves coding  nodes.
\item
$f$ preserves passing numbers at  coding nodes.
\end{enumerate}

Given a subtree $T\sse\bT_k$, let $G_T$ denote the graph represented by the coding nodes in $T$.
Notice that if $T\sse\mathbb{T}_k$ is strongly similar to $\mathbb{T}_k$, then $G_T$ is isomorphic to
 $\mathcal{H}_k$  as ordered graphs.

Essentially, a strong coding tree is a subtree of
$\mathbb{T}_k$ which is  strongly similar to
$\mathbb{T}_k$.
We say essentially, because there is one more important  consideration when working with  forbidden $k$-cliques.
Any  finite subtree of $\mathbb{T}_k$ which we build needs to be extendable within $\mathbb{T}_k$ to a subtree which is strongly similar to $\mathbb{T}_k$.
There are many finite subtrees
 of $\mathbb{T}_k$
for which this is not possible, because  each
remembers where it came from, coding edges with the original graph represented by the coding nodes in $\mathbb{T}_k$.

Take for example $k=3$.
If $A$ is a finite subtree of $\mathbb{T}_3$ and
 two  nodes $s,t\in A$ have passing number $1$ at some coding node $c_i$ in $\mathbb{T}_3$,
then whenever  $s$ is extended to some  coding node  $c_m$ in
$\mathbb{T}_3$,
then
any extension of $t$  to a coding node $c_n$ in $\bT_3$ with  length  greater than $|c_m|$  cannot have passing number $1$ at $c_m$, as that would code a triangle; no coding of a triangle occurs in
$\mathbb{T}_3$.
In such a situation,  there is no way to  extend $A$ to a subtree of $\mathbb{T}_3$ coding all of $\mathcal{H}_3$.

To prevent this, we make  a further requirement which is, roughly, as follows:
Fix $a\in [3,k]$.
We say that a level set $X\sse \bT_k$  with nodes of length $\ell_X$  {\em has a  pre-$a$-clique}  if  for some
$\mathcal{I} \sse[\om]^{a-2}$,
letting  $i_*=\max(\mathcal{I})$ and $\ell_*=
|c_{i_*}|$,
we have that
$\ell_*\le \ell_X$,
the set
$\{c_{i}:i\in \mathcal{I}\}$ codes an $(a-2)$-clique, and
each node in $X^+$
 has passing number $1$ at $c_{i}$, for each $i \in \mathcal{I}$.
The idea is that pre-$a$-cliques  code entanglements.
Essentially,
we say that a  subtree  $T$ of  $\bT_k$ has the {\em Witnessing Property}  if for each   pre-$a$-clique in $T$,
 $a\in [3,k]$,
there is a set of $(a-2)$ many coding nodes
$\{c_i:i\in\mathcal{I}\}$ as above, all of which are coding nodes in $T$.
 A tree $T\sse \bT_k$
is a member of $\mathcal{T}_k$ iff $T$
is strongly similar to $\mathbb{T}_k$ and
 has the Witnessing  Property (Lemma 5.15, \cite{DobrinenH_k19}).
Some examples of unwitnessed and witnessed pre-cliques are in Figures 7 -- 10.

Essentially, we define
  the space of {\em strong coding trees}, $\mathcal{T}_k$, to consist of those subtrees of $\mathbb{T}_k$ which are strongly similar  to $\mathbb{T}_k$ and have the Witnessing Property.
For  the details,  see \cite{DobrinenJML20} and \cite{DobrinenH_k19}.
This set of subtrees of $\mathcal{T}_k$ is the analogue of  strong trees appropriate to the Henson graph  $\mathcal{H}_k$



\begin{center}
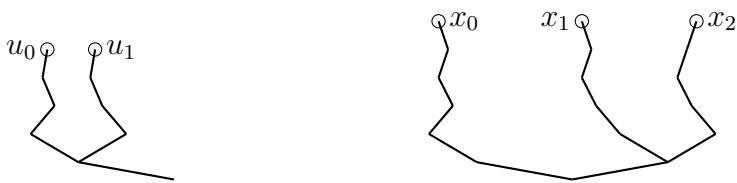
\begin{figure}
\begin{tikzpicture}[grow'=up,scale=.25]
\tikzstyle{level 1}=[sibling distance=4in]
\tikzstyle{level 2}=[sibling distance=2in]
\tikzstyle{level 3}=[sibling distance=1in]
\tikzstyle{level 4}=[sibling distance=0.5in]
\tikzstyle{level 5}=[sibling distance=0.2in]
\tikzstyle{level 6}=[sibling distance=0.2in]
\node {}
child{coordinate (t0) edge from parent[thick]
			child{coordinate (t00)
child{coordinate (t000)edge from parent[draw=none]
}
child{coordinate (t001)
child{coordinate (t0010)edge from parent[color=black]
child{coordinate (t00100)edge from parent[draw=none]}
child{coordinate (t00101)}
}
child{coordinate (t0011)edge from parent[draw=none]}}
}
			child{ coordinate(t01)
child{ coordinate(t010)
child{ coordinate(t0100)
child{ coordinate(t01000)edge from parent[draw=none]}
child{ coordinate(t01001)
}}
child{ coordinate(t0101)edge from parent[draw=none]
}
}
child{ coordinate(t011)edge from parent[draw=none]
}}}
child{ coordinate(t1)edge from parent[draw=none]
child{coordinate(t10)edge from parent[draw=none]
child{coordinate(t100)edge from parent[draw=none]
child{coordinate(t1000)edge from parent[draw=none]}
child{coordinate(t1001)edge from parent[draw=none]}
}
child{coordinate(t101)edge from parent[draw=none]}}
child{coordinate(t11)edge from parent[draw=none]}
			 };

\node[circle, draw, inner sep=0pt, minimum size=5pt] at (t00101)  {};

\node[circle, draw,inner sep=0pt, minimum size=5pt] at (t01001) {};
\node[left] at (t00101) {$u_0$};
\node[right] at (t01001) {$u_1$};
\end{tikzpicture}
\hspace{20pt}
\begin{tikzpicture}[grow'=up,scale=.25]
\tikzstyle{level 1}=[sibling distance=4in]
\tikzstyle{level 2}=[sibling distance=2in]
\tikzstyle{level 3}=[sibling distance=1in]
\tikzstyle{level 4}=[sibling distance=0.6in]
\tikzstyle{level 5}=[sibling distance=0.4in]
\node {}
child{coordinate (t0)edge from parent[thick]
			child{coordinate (t00)
child{coordinate (t000)edge from parent[draw=none]
child{coordinate (t0000)edge from parent[draw=none]}
child{coordinate(t0001)edge from parent[draw=none]}
}
child{ coordinate(t001)
child{ coordinate(t0010)
child{ coordinate(t00100)edge from parent[draw=none]}
child{ coordinate(t00101)
child{coordinate(t001010)edge from parent[color=black]}
child{coordinate(t001011)edge from parent[draw=none]} }}
child{ coordinate(t0011)edge from parent[draw=none]}}}
			child{ coordinate(t01)edge from parent[draw=none]
}}
		child{ coordinate(t1)edge from parent[thick]
			child{ coordinate(t10)
child{ coordinate(t100)
child{ coordinate(t1000)
child{ coordinate(t10000)edge from parent[draw=none]}
child{ coordinate(t10001)
child{coordinate(t100010)edge from parent[color=black]}
child{coordinate(t100011)edge from parent[draw=none]}
}}
child{ coordinate(t1001)edge from parent[draw=none]}}
child{ coordinate(t101)edge from parent[draw=none]}}
			child{  coordinate(t11)
child{ coordinate(t110)
child{ coordinate(t1100)
child{ coordinate(t11000)edge from parent[draw=none]}
child{ coordinate(t11001)
child{coordinate(t110010)edge from parent[draw=none]}
child{coordinate(t110011)edge from parent[color=black]}
}}
child{ coordinate(t1101)edge from parent[draw=none]}}
child{  coordinate(t111)edge from parent[draw=none]}}};

\node[circle, draw,inner sep=0pt, minimum size=5pt] at (t100010)  {};
\node[circle, draw, inner sep=0pt, minimum size=5pt] at (t001010)  {};
\node[circle, draw,inner sep=0pt, minimum size=5pt] at (t110011)  {};

\node[right] at (t001010) {$x_0$};
\node[left] at (t100010) {$x_1$};
\node[right] at (t110011) {$x_2$};

\end{tikzpicture}
\caption{Two examples of unwitnessed pre-$3$-cliques}
\end{figure}
\end{center}


\begin{center}
\begin{figure}
\begin{tikzpicture}[grow'=up,scale=.25]
\tikzstyle{level 1}=[sibling distance=4in]
\tikzstyle{level 2}=[sibling distance=2in]
\tikzstyle{level 3}=[sibling distance=1in]
\tikzstyle{level 4}=[sibling distance=0.5in]
\tikzstyle{level 5}=[sibling distance=0.2in]
\tikzstyle{level 6}=[sibling distance=0.2in]
\node {}
child{coordinate (t0) edge from parent[thick]
			child{coordinate (t00)
child{coordinate (t000)edge from parent[draw=none]
}
child{coordinate (t001)
child{coordinate (t0010)edge from parent[color=black]
child{coordinate (t00100)edge from parent[draw=none]}
child{coordinate (t00101)}
}
child{coordinate (t0011)edge from parent[draw=none]}}
}
			child{ coordinate(t01)
child{ coordinate(t010)
child{ coordinate(t0100)
child{ coordinate(t01000)edge from parent[draw=none]}
child{ coordinate(t01001)
}}
child{ coordinate(t0101)edge from parent[draw=none]
}
}
child{ coordinate(t011)edge from parent[draw=none]
}}}
child{ coordinate(t1)edge from parent[thick]
child{coordinate(t10)
child{coordinate(t100)
child{coordinate(t1000)edge from parent[draw=none]}
child{coordinate(t1001)}}
child{coordinate(t101)edge from parent[draw=none]}}
child{coordinate(t11)edge from parent[draw=none]}
			 };

\node[circle, draw, inner sep=0pt, minimum size=5pt] at (t00101)  {};
\node[circle, draw,inner sep=0pt, minimum size=5pt] at (t01001) {};
\node[circle, fill=black,inner sep=0pt, minimum size=5pt] at (t1001)  {};

\node[left] at (t00101) {$u_0$};
\node[right] at (t01001) {$u_1$};
\node[right] at (t1001) {$c_n$};

\end{tikzpicture}
\hspace{20pt}
\begin{tikzpicture}[grow'=up,scale=.25]
\tikzstyle{level 1}=[sibling distance=4in]
\tikzstyle{level 2}=[sibling distance=2in]
\tikzstyle{level 3}=[sibling distance=1in]
\tikzstyle{level 4}=[sibling distance=0.6in]
\tikzstyle{level 5}=[sibling distance=0.4in]
\node {}
child{coordinate (t0)edge from parent[thick]
			child{coordinate (t00)
child{coordinate (t000)
child{coordinate (t0000)edge from parent[draw=none]}
child{coordinate(t0001)}
}
child{ coordinate(t001)
child{ coordinate(t0010)
child{ coordinate(t00100)edge from parent[draw=none]}
child{ coordinate(t00101)
child{coordinate(t001010)edge from parent[color=black]}
child{coordinate(t001011)edge from parent[draw=none]} }}
child{ coordinate(t0011)edge from parent[draw=none]}}}
			child{ coordinate(t01)edge from parent[draw=none]
}}
		child{ coordinate(t1)edge from parent[thick]
			child{ coordinate(t10)
child{ coordinate(t100)
child{ coordinate(t1000)
child{ coordinate(t10000)edge from parent[draw=none]}
child{ coordinate(t10001)
child{coordinate(t100010)edge from parent[color=black]}
child{coordinate(t100011)edge from parent[draw=none]}
}}
child{ coordinate(t1001)edge from parent[draw=none]}}
child{ coordinate(t101)edge from parent[draw=none]}}
			child{  coordinate(t11)
child{ coordinate(t110)
child{ coordinate(t1100)
child{ coordinate(t11000)edge from parent[draw=none]}
child{ coordinate(t11001)
child{coordinate(t110010)edge from parent[draw=none]}
child{coordinate(t110011)edge from parent[color=black]}
}}
child{ coordinate(t1101)edge from parent[draw=none]}}
child{  coordinate(t111)edge from parent[draw=none]}}};

\node[circle, draw,inner sep=0pt, minimum size=5pt] at (t100010)  {};
\node[circle, draw, inner sep=0pt, minimum size=5pt] at (t001010)  {};
\node[circle, draw,inner sep=0pt, minimum size=5pt] at (t110011)  {};
\node[circle, fill=black,inner sep=0pt, minimum size=5pt] at (t0001) {};

\node[right] at (t001010) {$x_0$};
\node[left] at (t100010) {$x_1$};
\node[right] at (t110011) {$x_2$};
\node[left] at (t0001) {$c_n$};

\end{tikzpicture}
\caption{Two examples of witnessed pre-$3$-cliques}
\end{figure}
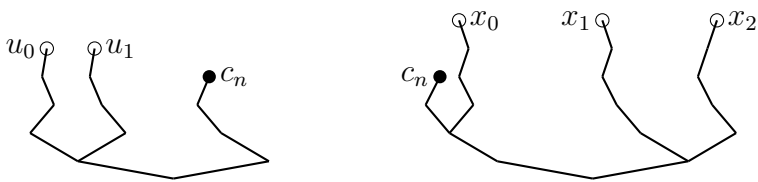
\end{center}



\begin{center}
\begin{figure}
\begin{tikzpicture}[grow'=up,scale=.2]
\tikzstyle{level 1}=[sibling distance=4in]
\tikzstyle{level 2}=[sibling distance=2in]
\tikzstyle{level 3}=[sibling distance=1in]
\tikzstyle{level 4}=[sibling distance=0.6in]
\tikzstyle{level 5}=[sibling distance=0.4in]
\node {}
child{coordinate (t0)edge from parent[thick]
			child{coordinate (t00)
child{coordinate (t000)edge from parent[draw=none]}
child{ coordinate(t001)
child{ coordinate(t0010)
child{ coordinate(t00100)edge from parent[draw=none]}
child{ coordinate(t00101)
child{coordinate(t001010)edge from parent[color=black]
child{coordinate(t0010100)edge from parent[draw=none]}
child{coordinate(t0010101)}}
child{coordinate(t001011)edge from parent[draw=none]} }}
child{ coordinate(t0011)edge from parent[draw=none]}}}
			child{ coordinate(t01)edge from parent[draw=none]
}}
		child{ coordinate(t1)edge from parent[thick]
			child{ coordinate(t10)
child{ coordinate(t100)
child{ coordinate(t1000)
child{ coordinate(t10000)edge from parent[draw=none]}
child{ coordinate(t10001)
child{coordinate(t100010)edge from parent[color=black]
child{coordinate(t1000100)edge from parent[draw=none]}
child{coordinate(t1000101)}}
child{coordinate(t100011)edge from parent[draw=none]}
}}
child{ coordinate(t1001)edge from parent[draw=none]}}
child{ coordinate(t101)edge from parent[draw=none]}}
			child{  coordinate(t11)edge from parent[draw=none]
}};

\node[circle, draw,inner sep=0pt, minimum size=5pt] at (t1000101)  {};
\node[circle, draw, inner sep=0pt, minimum size=5pt] at (t0010101)  {};

\node[right] at (t0010101) {$y_0$};
\node[left] at (t1000101) {$y_1$};

\end{tikzpicture}
\hspace{20pt}
\begin{tikzpicture}[grow'=up,scale=.2]
\tikzstyle{level 1}=[sibling distance=4in]
\tikzstyle{level 2}=[sibling distance=2in]
\tikzstyle{level 3}=[sibling distance=1in]
\tikzstyle{level 4}=[sibling distance=0.6in]
\tikzstyle{level 5}=[sibling distance=0.4in]
\node {}
child{coordinate (t0)edge from parent[thick]
			child{coordinate (t00)
child{coordinate (t000)edge from parent[draw=none]}
child{ coordinate(t001)
child{ coordinate(t0010)
child{ coordinate(t00100)edge from parent[draw=none]}
child{ coordinate(t00101)
child{coordinate(t001010)edge from parent[color=black]
child{coordinate(t0010100)edge from parent[draw=none]}
child{coordinate(t0010101)}}
child{coordinate(t001011)edge from parent[draw=none]} }}
child{ coordinate(t0011)edge from parent[draw=none]}}}
			child{ coordinate(t01)edge from parent[draw=none]
}}
		child{ coordinate(t1)edge from parent[thick]
			child{ coordinate(t10)
child{ coordinate(t100)
child{ coordinate(t1000)
child{ coordinate(t10000)edge from parent[draw=none]}
child{ coordinate(t10001)
child{coordinate(t100010)edge from parent[color=black]
child{coordinate(t1000100)edge from parent[draw=none]}
child{coordinate(t1000101)}}
child{coordinate(t100011)edge from parent[draw=none]}
}}
child{ coordinate(t1001)edge from parent[draw=none]}}
child{ coordinate(t101)edge from parent[draw=none]}}
			child{  coordinate(t11)
child{ coordinate(t110)
child{ coordinate(t1100)
child{ coordinate(t11000)edge from parent[draw=none]}
child{ coordinate(t11001)
child{coordinate(t110010)edge from parent[color=black]
child{coordinate(t1100100)edge from parent[draw=none]}
child{coordinate(t1100101)}
}
child{coordinate(t110011)edge from parent[draw=none]}
}}
child{ coordinate(t1101)edge from parent[draw=none]}}
child{  coordinate(t111)edge from parent[draw=none]}}};

\node[circle, draw,inner sep=0pt, minimum size=5pt] at (t1000101)  {};
\node[circle, draw, inner sep=0pt, minimum size=5pt] at (t0010101)  {};
\node[circle, draw,inner sep=0pt, minimum size=5pt] at (t1100101)  {};

\node[right] at (t0010101) {$z_0$};
\node[left] at (t1000101) {$z_1$};
\node[right] at (t1100101) {$z_2$};

\end{tikzpicture}
\caption{Two examples of unwitnessed pre-$4$-cliques}
\end{figure}
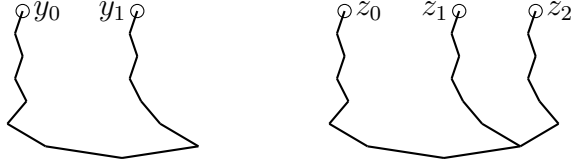
\end{center}


\begin{center}
\begin{figure}
\begin{tikzpicture}[grow'=up,scale=.2]
\tikzstyle{level 1}=[sibling distance=4in]
\tikzstyle{level 2}=[sibling distance=2in]
\tikzstyle{level 3}=[sibling distance=1in]
\tikzstyle{level 4}=[sibling distance=0.6in]
\tikzstyle{level 5}=[sibling distance=0.4in]
\node {}
child{coordinate (t0)edge from parent[thick]
			child{coordinate (t00)
child{coordinate (t000)
child{coordinate (t0000)
child{coordinate(t00000)edge from parent[draw=none]}
child{coordinate(t00001)
child{coordinate(t000010)edge from parent[draw=none]}
child{coordinate(t000011)}}}
child{coordinate(t0001)}
}
child{ coordinate(t001)
child{ coordinate(t0010)
child{ coordinate(t00100)edge from parent[draw=none]}
child{ coordinate(t00101)
child{coordinate(t001010)edge from parent[color=black]
child{coordinate(t0010100)edge from parent[draw=none]}
child{coordinate(t0010101)}}
child{coordinate(t001011)edge from parent[draw=none]} }}
child{ coordinate(t0011)edge from parent[draw=none]}}}
			child{ coordinate(t01)edge from parent[draw=none]
}}
		child{ coordinate(t1)edge from parent[thick]
			child{ coordinate(t10)
child{ coordinate(t100)
child{ coordinate(t1000)
child{ coordinate(t10000)edge from parent[draw=none]}
child{ coordinate(t10001)
child{coordinate(t100010)edge from parent[color=black]
child{coordinate(t1000100)edge from parent[draw=none]}
child{coordinate(t1000101)}}
child{coordinate(t100011)edge from parent[draw=none]}
}}
child{ coordinate(t1001)edge from parent[draw=none]}}
child{ coordinate(t101)edge from parent[draw=none]}}
			child{  coordinate(t11)edge from parent[draw=none]
}};

\node[circle, draw,inner sep=0pt, minimum size=5pt] at (t1000101)  {};
\node[circle, draw, inner sep=0pt, minimum size=5pt] at (t0010101)  {};
\node[circle, fill=black,inner sep=0pt, minimum size=5pt] at (t0001) {};
\node[circle, fill=black,inner sep=0pt, minimum size=5pt] at (t000011) {};

\node[right] at (t0010101) {$y_0$};
\node[left] at (t1000101) {$y_1$};
\node[left] at (t0001) {$c_m$};
\node[left] at (t000011) {$c_n$};

\end{tikzpicture}
\hspace{20pt}
\begin{tikzpicture}[grow'=up,scale=.2]
\tikzstyle{level 1}=[sibling distance=4in]
\tikzstyle{level 2}=[sibling distance=2in]
\tikzstyle{level 3}=[sibling distance=1in]
\tikzstyle{level 4}=[sibling distance=0.6in]
\tikzstyle{level 5}=[sibling distance=0.4in]
\node {}
child{coordinate (t0)edge from parent[thick]
			child{coordinate (t00)
child{coordinate (t000)
child{coordinate (t0000)
child{coordinate(t00000)edge from parent[draw=none]}
child{coordinate(t00001)
child{coordinate(t000010)edge from parent[draw=none]}
child{coordinate(t000011)}}}
child{coordinate(t0001)}
}
child{ coordinate(t001)
child{ coordinate(t0010)
child{ coordinate(t00100)edge from parent[draw=none]}
child{ coordinate(t00101)
child{coordinate(t001010)edge from parent[color=black]
child{coordinate(t0010100)edge from parent[draw=none]}
child{coordinate(t0010101)}}
child{coordinate(t001011)edge from parent[draw=none]} }}
child{ coordinate(t0011)edge from parent[draw=none]}}}
			child{ coordinate(t01)edge from parent[draw=none]
}}
		child{ coordinate(t1)edge from parent[thick]
			child{ coordinate(t10)
child{ coordinate(t100)
child{ coordinate(t1000)
child{ coordinate(t10000)edge from parent[draw=none]}
child{ coordinate(t10001)
child{coordinate(t100010)edge from parent[color=black]
child{coordinate(t1000100)edge from parent[draw=none]}
child{coordinate(t1000101)}}
child{coordinate(t100011)edge from parent[draw=none]}
}}
child{ coordinate(t1001)edge from parent[draw=none]}}
child{ coordinate(t101)edge from parent[draw=none]}}
			child{  coordinate(t11)
child{ coordinate(t110)
child{ coordinate(t1100)
child{ coordinate(t11000)edge from parent[draw=none]}
child{ coordinate(t11001)
child{coordinate(t110010)edge from parent[color=black]
child{coordinate(t1100100)edge from parent[draw=none]}
child{coordinate(t1100101)}
}
child{coordinate(t110011)edge from parent[draw=none]}
}}
child{ coordinate(t1101)edge from parent[draw=none]}}
child{  coordinate(t111)edge from parent[draw=none]}}};

\node[circle, draw,inner sep=0pt, minimum size=5pt] at (t1000101)  {};
\node[circle, draw, inner sep=0pt, minimum size=5pt] at (t0010101)  {};
\node[circle, draw,inner sep=0pt, minimum size=5pt] at (t1100101)  {};
\node[circle, fill=black,inner sep=0pt, minimum size=5pt] at (t0001) {};
\node[circle, fill=black,inner sep=0pt, minimum size=5pt] at (t000011) {};

\node[right] at (t0010101) {$z_0$};
\node[left] at (t1000101) {$z_1$};
\node[right] at (t1100101) {$z_2$};
\node[left] at (t0001) {$c_m$};
\node[left] at (t000011) {$c_n$};

\end{tikzpicture}
\caption{Two examples of witnessed pre-$4$-cliques}
\end{figure}
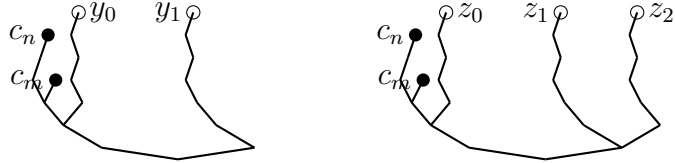
\end{center}


\subsection{A Ramsey theorem for strictly similar antichains}

The upper bounds for big Ramsey degrees of the Henson graphs basically  come from the fact that there are only finitely many ways to code a given finite graph within a strong coding tree.
Analogously to the case for the Rado graph, given any strong coding tree $T\in\mathcal{T}_k$, there is an antichain of coding nodes in $T$ which code the Henson graph $\mathcal{H}_k$.
Thus, one can restrict attention to colorings of antichains  representing a given finite $k$-clique-free graph, say $G$.
The relevant structural properties are those of strong similarity and new pre-cliques and their placement  within the antichain.
This is the idea behind {\em strict similarity}.  For the precise definition, the reader is referred to \cite{DobrinenJML20} and \cite{DobrinenH_k19}.

\begin{thm}[\cite{DobrinenJML20},\cite{DobrinenH_k19}]\label{thm.SWPAntichain}
Let $Z$ be a finite antichain of coding nodes in  a strong $\mathcal{H}_k$-coding tree $T\in\mathcal{T}_k$,
and suppose  $h$  colors of all antichains  $T$ which are strictly similar to $Z$ into finitely many colors.
Then there is an  strong $\mathcal{H}_k$-coding tree $S\le T$ such that
all subsets of $S$  strictly similar to $Z$ have the same $h$ color.
\end{thm}

The proof uses forcing to obtain a ZFC result, but not in the usual manner using absoluteness.
Recall the Halpern-\Lauchli\ Theorem \ref{thm.HL}.
Harrington gave an insightful  proof which uses Cohen forcing in the following way.
Suppose we have $d$ many infinite strong trees $T_i$ and a coloring of the  product of their level sets as in the statement of Theorem \ref{thm.HL}.
Let $\kappa$ be large enough that $\kappa\ra(\aleph_1)^{2d}_{\aleph_0}$ and take $\mathbb{P}$ to be $\kappa$-Cohen forcing which adds $\kappa$ new branches to each of the $d$ many trees.
Harrington gave an argument  guaranteeing  that there are nodes $t_i^*\in T_i$, for each $i<d$, of the same length which have some color $\varepsilon\in r$;
moreover, given any level sets extending  these nodes, there are further extensions to level sets so that each member in their product has the same color $\varepsilon$.
The forcing is used to find these finite sets successively in $\om$ many steps;  the generic filter is never actually used - one never actually passes to a generic extension.
Instead, the forcing language and basic facts about forcing guarantee that certain finite level sets exist, and their finiteness guarantees that they are actually in the ground model.
So, this is very much not a constructive proof, but at the same time, it is a ZFC proof where one constructs each level of the subtree separately.

The author extended this idea  to the context strong coding trees.
However, since  strong coding trees have two types of nodes, coding and non-coding, new forcings which are not equivalent to $\kappa$-Cohen forcing had to be introduced.
The general approach uses Harrington's ideas as a starting point, but the implementation is much more involved.
Another  new element in this setting is that envelopes will be antichains of coding nodes which have a strong version of the Witnessing Property.
This is quite different from envelopes for the rationals or the Rado graph being finite strong trees.

Theorem \ref{thm.SWPAntichain} is applied as follows to prove the finite big Ramsey degrees.
Given a finite $k$-clique-free graph $G$,
there are only finitely many strict similarity types of antichains of coding nodes representing $G$.
The number of such strict similarity types is the upper bound for the big Ramsey degree of $G$ in $\mathcal{H}_k$.

\begin{thm}[\cite{DobrinenJML20},\cite{DobrinenH_k19}]\label{thm.main}
Suppose  $k\ge 3$ and $G$ is a finite $k$-clique-free graph.
Let $h$ color all copies of $G$ in $\mathcal{H}_k$ into finitely many colors.
Then there is a subgraph of $\mathcal{H}_k$ which is isomorphic to $\mathcal{H}_k$ in which the copies of $G$ take on no more colors than the number of strict similarity types of antichains in $\bT_k$ coding $G$.
\end{thm}


\subsection{Infinite dimensional Ramsey theory of the Rado graph}

We now mention a recent result of the author on infinite dimensional Ramsey theory of the Rado graph.
In  Problem 11.2 of  \cite{Kechris/Pestov/Todorcevic05}, Kechris, Pestov, and Todorcevic
 ask  for  the topological dynamics analogue
 of a    corresponding   infinite  Ramsey-theoretic   result for several \Fraisse\ structures, in particular,
  the rationals, the Rado graph, and  the Henson graphs.
By an infinite Ramsey-theoretic result,
they mean
 a result of the form
\begin{equation}\label{eq.Prob11.2}
\mathbb{F}\ra_* (\mathbb{F})^{\mathbb{F}}_{l,t},
\end{equation}
where  equation (\ref{eq.Prob11.2})
reads:
``For each  partition of
 ${\mathbb{F}\choose\mathbb{F}}$ into $l$ many  definable   subsets,
there is an $\mathbf{F}\in {\mathbb{F}\choose\mathbb{F}}$
such that ${\mathbf{F}\choose\mathbb{F}}$ is  contained in
no more than $t$ of the pieces of the partition.''
Here,
one  assumes a  natural topology on  ${\mathbb{F}\choose\mathbb{F}}$ and
{\em definable} refers to
any reasonable class of sets definable relative to  the topology, for instance,
 open, Borel, analytic,  or property of Baire.
A sub-question   implicit in   Problem 11.2 in \cite{Kechris/Pestov/Todorcevic05} is the following broader version of Question \ref{question.2.9}:

\begin{question}\label{q.KPT}
For which ultrahomogeneous structures $\mathbb{F}$
does it hold that
\begin{equation}\label{eq.qKPT}
\mathbb{F}\,  \ra_* \, (\mathbb{F} )^{\mathbb{F}}_{l,t},
\end{equation}
for some positive integer $t$?
\end{question}

The natural topology  to
give such a space is the
 one induced by ordering  the universe $F$ of $\mathbb{F}$  in order-type $\om$, and viewing
 ${\mathbb{F}\choose\mathbb{F}}$ as a subspace of the product space $2^{F}$  with the Tychonoff topology.
 Kechris, Pestov, and Todorcevic pointed out  that  very little is known about Question \ref{q.KPT}.

 In \cite{DobrinenRado19}, the author set out to answer this question for the Rado graph.
Since
 the big Ramsey degrees for copies of a finite graph inside  the Rado graph grow  without bound as the number of vertices in the finite graph  whose copies are being colored  grows  implies that
any positive  answer to Question \ref{q.KPT} for the Rado graph
 must restrict to  a collection of Rado graphs all of whose vertices are ordered in the same order.
Furthermore,
 it is necessary that all copies of the Rado graph being colored have the same strong similarity type.
 Otherwise, one
 may use strong similarity types to make a coloring of the copies of the Rado graph to show that there is no bound  $t$ of the sort in   (\ref{eq.Prob11.2}), where $\mathbb{F}$ is the Rado graph.

In \cite{DobrinenRado19},
the author  answered  Question \ref{q.KPT} for
a  collection of  Rado graphs, each of which has  the same strong similarity type.
While the Rado graph can be represented by nodes in strong trees of the kind in the Milliken Theorem \ref{thm.Milliken},
that theorem by itself does not answer this question, as it is unclear in a tree without coding nodes how the strong subtrees should be thought of as coding subcopies of a given Rado graph.
In order to make sure that the representations of the subgraphs were concrete, it turned out to be
useful to work with trees with coding nodes, even though the Rado graph has no forbidden subgraphs.

Let $\mathbb{R}=(R,E)$ denote the Rado graph with vertices ordered as $\lgl v_n:n<\om\rgl$  represented by the coding nodes $\lgl c_n:n<\om\rgl$ in the tree $\bT_{\mathbb{R}}$ in Figure 11.
Let
$\mathcal{T}_{\mathbb{R}}$ denote the  collection of all subtrees
$T\sse \mathbb{T}_{\mathbb{R}}$  strongly similar to $\mathbb{T}_{\mathbb{R}}$.
Given $T\in \mathcal{T}_{\mathbb{R}}$, let $G_T$ denote the ordered graph represented by the coding nodes of $T$, noting that each $G_T$ is a subcopy of the Rado graph, ordered in the same way as $\mathbb{R}$.
Let $\mathcal{R}$ be the collection of all $G_T$, where $T\in \mathcal{T}_{\mathbb{R}}$.
 The topology on $\mathcal{R}$ is the topology inherited from the Tychonoff topology on
 $2^{R}$.

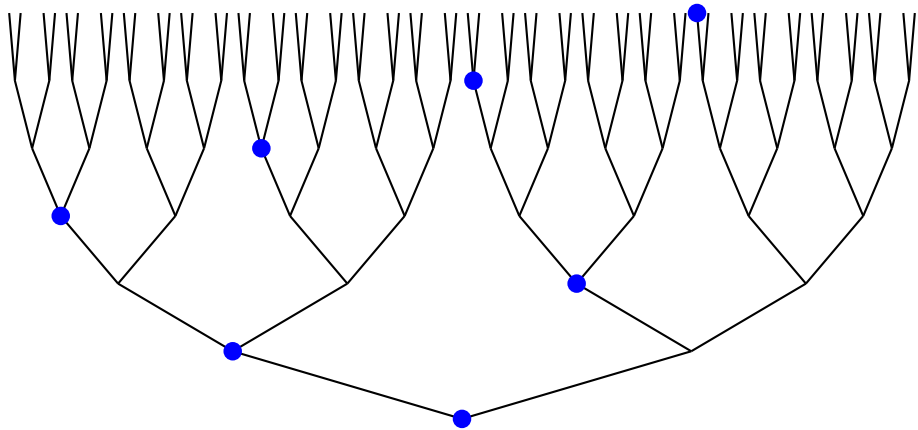
\begin{figure}
\begin{tikzpicture}[grow'=up,scale=.6]
\tikzstyle{level 1}=[sibling distance=4in]
\tikzstyle{level 2}=[sibling distance=2in]
\tikzstyle{level 3}=[sibling distance=1in]
\tikzstyle{level 4}=[sibling distance=0.5in]
\tikzstyle{level 5}=[sibling distance=0.3in]
\tikzstyle{level 6}=[sibling distance=0.1in]
\node {} coordinate (t9)
child{coordinate (t0) edge from parent[thick]
			child{coordinate (t00)
child{coordinate (t000)
child {coordinate(t0000)edge from parent[color=black]
child{coordinate(t00000)
child
child}
child{coordinate(t00001)
child
child}}
child {coordinate(t0001)
child {coordinate(t00010)
child
child }
child{coordinate(t00011)
child
child}}}
child{ coordinate(t001)
child{ coordinate(t0010)
child{ coordinate(t00100)
child
child }
child{ coordinate(t00101)
child
child }}
child{ coordinate(t0011) edge from parent[color=black]
child{ coordinate(t00110)
child
child
}
child{ coordinate(t00111)
child
child}}}}
			child{ coordinate(t01)
child{ coordinate(t010)
child{ coordinate(t0100)
child{ coordinate(t01000)
child
child}
child{ coordinate(t01001)
child
child }}
child{ coordinate(t0101)edge from parent[color=black]
child{ coordinate(t01010)
child
child}
child{ coordinate(t01011)
child
child}}}
child{ coordinate(t011)
child{ coordinate(t0110)
child{ coordinate(t01100)
child
child }
child{ coordinate(t01101)
child
child }}
child{ coordinate(t0111) edge from parent[color=black]
child { coordinate(t01110)
child
child}
child{ coordinate(t01111)
child
child}}}}}
		child{ coordinate(t1) edge from parent[thick, color=black]
			child{ coordinate(t10)
child{ coordinate(t100)
child{ coordinate(t1000) edge from parent[color=black]
child{ coordinate(t10000)
child
child}
child{ coordinate(t10001)
child
child}}
child{ coordinate(t1001)
child{ coordinate(t10010)
child
child}
child{ coordinate(t10011)
child
child}}}
child{ coordinate(t101)
child{ coordinate(t1010) edge from parent[color=black]
child{ coordinate(t10100)
child
child}
child{ coordinate(t10101)
child
child}}
child{ coordinate(t1011)
child{ coordinate(t10110)
child
child}
child{ coordinate(t10111)
child
child}}}}
			child{  coordinate(t11)  edge from parent[color=black]
child{ coordinate(t110)
child{ coordinate(t1100)
child{ coordinate(t11000)
child{coordinate(t110000)}
child}
child{ coordinate(t11001)
child
child}}
child{ coordinate(t1101)
child{ coordinate(t11010)
child
child}
child{ coordinate(t11011)
child
child}}}
child{  coordinate(t111)
child{  coordinate(t1110)
child{  coordinate(t11100)
child
child}
child{  coordinate(t11101)
child
child}}
child{  coordinate(t1111)
child{  coordinate(t11110)
child
child}
child{  coordinate(t11111)
child
child}}}} };

\node[circle, fill=blue,inner sep=0pt, minimum size=7pt] at (t0) {};
\node[circle, fill=blue,inner sep=0pt, minimum size=7pt] at (t10) {};
\node[circle, fill=blue,inner sep=0pt, minimum size=7pt] at (t9) {};
\node[circle, fill=blue,inner sep=0pt, minimum size=7pt] at (t000) {};
\node[circle, fill=blue,inner sep=0pt, minimum size=7pt] at (t0100) {};
\node[circle, fill=blue,inner sep=0pt, minimum size=7pt] at (t10000) {};
\node[circle, fill=blue,inner sep=0pt, minimum size=7pt] at (t110000) {};

\end{tikzpicture}
\caption{The Strong Rado Coding Tree $\bT_{\mathbb{R}}$}
\end{figure}

\begin{thm}[\cite{DobrinenRado19}]\label{thm.DRado}
If $\mathcal{X}\sse\mathcal{R}$ is Borel,
then for each graph $\mathbb{G}\in\mathcal{R}$,
either all members of $\mathcal{R}$ contained in $\mathbb{G}$ are members of $\mathcal{X}$,
or else
 no member of $\mathcal{R}$ contained in $\mathbb{G}$ is a  member of $\mathcal{X}$.
\end{thm}

A few remarks about this theorem are in order.
First, we point out that Theorem \ref{thm.DRado} is the analogue of the Galvin-Prikry Theorem \ref{thm.GP} for colorings of those copies of the Rado graph which have induced trees strongly similar to $\bT_{\mathbb{R}}$.
Second, the methods of proof use forcing in a similar yet simpler way than that in \cite{DobrinenJML20}.
However, these methods were not conducive to obtaining an analogue of the Ellentuck Theorem  \ref{thm.Ellentuck}.
So it remains open whether or not, with the natural Ellentuck-like topology on $\mathcal{R}$,
the subsets of $\mathcal{R}$ with the property of Baire have the Ramsey property.

Third, the strong coding tree $\mathbb{T}_{\mathbb{R}}$ is  the simplest kind of tree coding $\mathbb{R}$.
One could, however, fix any  perfect tree (skew or not), say $\mathbb{T}'$,  with coding nodes  which are dense and code the Rado graph, and the same methods would produce the same result for those subcopies of the Rado graph induced by the coding nodes of  subtrees of $\mathbb{T}'$ which
are strongly similar to $\mathbb{T}'$.
Lastly, the methods should produce  a similar theorem for any of the binary relational structures in \cite{Sauer06} and \cite{Laflamme/Sauer/Vuksanovic06}.


\section{Future directions}\label{sec.fd}

As seen in the previous section, trees with coding nodes have provided means for solving problems on big Ramsey degrees for Henson graphs as well as infinite dimensional Ramsey theory for the Rado graph, which has no forbidden subgraphs.
We envision several directions in which this idea can be developed to solve questions on Ramsey theory of infinite structures.

It is  likely that trees with coding nodes will be useful  in proving  finite big Ramsey degrees for  all  \Fraisse\ limits of \Fraisse\ classes with finitely many binary relations
and  a finite constraint set
 satisfying the Ramsey property.
Indeed, it may be enough to require finite small Ramsey degrees.
Let $\mathcal{K}$ be a  \Fraisse\ class in 
 a finite  relational language $\mathcal{L}$, where all relations are finitary.
 Let us call a  set $\mathbf{F}$ of finite structures in $\mathcal{L}$ a {\em  constraint set}  for $\mathcal{K}$ if
 for any finite structure $\bsA$ in the language $\mathcal{L}$,
  $\bsA\in \mathcal{K}$ if and only if no induced substructure of $\bsA$ is in $\mathbf{F}$.
In \cite{Sauer03},
Sauer  gave a detailed characterization of the big Ramsey degrees of vertices in transitive free amalgamation structures in finite binary languages
 and provided some explicit examples of ultrahomogeneous directed graphs where the vertices do not have finite big Ramsey degrees.
By the \Nesetril-\Rodl\ Theorem, one can consider the linearly ordered versions of these classes to obtain Ramsey classes whose vertices do not have finite big Ramsey degrees.

At the  conference, {\em Unifying Themes in Ramsey Theory}  in Banff, 2018,
Sauer suggested
 trying  to move the forcing proofs from strong coding trees to a direct approach by forcing directly on the structures.
Recall that  the nodes in a strong coding tree represent a finite partial $1$-type over a finite substructure; that is, the trees are really just a means for visualizing or making a structure out of the types.
 One important task is to interpret the forcings  in the work related in Section 4    back into the graph setting in a way   that points  to the natural     analogues of strong coding trees for structures with relations of any
finite arity.
Such an approach can hopefully lead to proving  the following conjecture.

\begin{conjecture}
Let $\mathcal{K}$ be a relational  \Fraisse\ class 
 with  finitely many relations and  a finite constraint set.
Suppose that $\mathcal{K}$ satisfies 
the Ramsey property (or has finite small Ramsey degrees).
Then  the \Fraisse\ limit of $\mathcal{K}$ has finite big Ramsey degrees. 
\end{conjecture}

The methods for the infinite dimensional Ramsey theory of the Rado graph are a simplified version of those developed for the big Ramsey degrees of the Henson graphs.
I seems to be just a matter of double checking to see that the work in \cite{DobrinenH_k19} will induce
 infinite dimensional Ramsey theory for Borel sets   in  any class of   $k$-clique-free Henson graphs.
It will be interesting to see if some other methods can produce analogues of the Ellentuck Theorem for these graphs, or whether there is some essential property of these graphs which prevent this.
It seems likely that whatever structures have finite big Ramsey degrees will also have infinite dimensional Ramsey theorems.

\begin{conjecture}\label{conj2}
 Let $\mathcal{K}$ be a relational  \Fraisse\ class 
 with  finitely many relations and  a finite constraint set.
Suppose that $\mathcal{K}$ satisfies 
the Ramsey property (or has finite small Ramsey degrees).
 Let $\mathbb{K}$ be the \Fraisse\ limit of $\mathcal{K}$ with some linear order of its universe in order-type $\om$.
 Then there is some notion of strong similarity type so that
 the collection  of all subcopies  of $\mathbb{K}$ having the same strong similarity type has the property that all Borel subsets are Ramsey.
\end{conjecture}


\bibliographystyle{amsplain}
\bibliography{references}

\end{document}